\documentclass{amsart}
\usepackage{microtype}
\usepackage{standalone}
\usepackage{amssymb,mathtools,mathbbol}
\usepackage{tikz-cd}
\usepackage{csquotes}
\usepackage{autobreak}

\usepackage[
	backend=biber,
	style=alphabetic,
	backref=true,
	url=false,
	doi=false,
	isbn=false,
	eprint=false]{biblatex}

\renewbibmacro{in:}{} 
\AtEveryBibitem{\clearfield{pages}} 

\DeclareFieldFormat{title}{\myhref{\mkbibemph{#1}}}
\DeclareFieldFormat
[article,inbook,incollection,inproceedings,patent,thesis,unpublished]
{title}{\myhref{\mkbibquote{#1\isdot}}}

\newcommand{\doiorurl}{%
	\iffieldundef{url}
	{\iffieldundef{eprint}
		{}
		{http://arxiv.org/abs/\strfield{eprint}}}
	{\strfield{url}}%
}

\newcommand{\myhref}[1]{%
	\ifboolexpr{%
		test {\ifhyperref}
		and
		not test {\iftoggle{bbx:eprint}}
		and
		not test {\iftoggle{bbx:url}}
	}
	{\href{\doiorurl}{#1}}
	{#1}%
}

\usepackage[
	bookmarks=true,
	linktocpage=true,
	bookmarksnumbered=true,
	breaklinks=true,
	pdfstartview=FitH,
	hyperfigures=false,
	plainpages=false,
	naturalnames=true,
	colorlinks=true,
	pagebackref=false,
	pdfpagelabels]{hyperref}

\hypersetup{
	colorlinks,
	citecolor=blue,
	filecolor=blue,
	linkcolor=blue,
	urlcolor=blue
}

\usepackage[capitalize, noabbrev]{cleveref}
\crefname{subsection}{\S\!}{subsections}

\makeatletter
\@namedef{subjclassname@2020}{%
	\textup{2020} Mathematics Subject Classification}
\makeatother

\newtheorem{theorem}{Theorem}
\newtheorem*{theorem*}{Theorem}
\newtheorem{proposition}[theorem]{Proposition}
\newtheorem*{proposition*}{Proposition}
\newtheorem{lemma}[theorem]{Lemma}
\newtheorem*{lemma*}{Lemma}
\newtheorem{corollary}[theorem]{Corollary}
\newtheorem*{corollary*}{Corollary}

\theoremstyle{definition}
\newtheorem{definition}[theorem]{Definition}
\newtheorem*{definition*}{Definition}

\newtheorem*{remark*}{Remark}

\newtheorem*{example*}{Example}

\newtheorem*{construction*}{Construction}

\newtheorem*{convention*}{Convention}

\newtheorem*{terminology*}{Terminology}

\newtheorem*{notation*}{Notation}

\newtheorem*{question*}{Question}

\DeclareMathOperator{\bd}{\partial}

\newcommand{\ot}{\otimes}
\DeclareMathOperator{\EZ}{EZ}

\newcommand{\N}{\mathbb{N}}
\newcommand{\Z}{\mathbb{Z}}

\renewcommand{\k}{\Bbbk}
\newcommand{\sym}{\mathbb{S}}

\newcommand{\gsimplex}{\mathbb{\Delta}}
\newcommand{\gcube}{\mathbb{I}}

\newcommand{\Cat}{\mathsf{Cat}}
\newcommand{\Fun}{\mathsf{Fun}}
\newcommand{\Set}{\mathsf{Set}}
\newcommand{\Top}{\mathsf{Top}}
\newcommand{\CW}{\mathsf{CW}}
\newcommand{\Ch}{\mathsf{Ch}}
\newcommand{\simplex}{\triangle}
\newcommand{\sSet}{\mathsf{sSet}}
\newcommand{\cube}{\square}
\newcommand{\cSet}{\mathsf{cSet}}

\newcommand{\coAlg}{\mathsf{coAlg}}
\newcommand{\biAlg}{\mathsf{biAlg}}

\newcommand{\symMod}{\mathsf{Mod}_{\sym}}

\newcommand{\operads}{\mathsf{Oper}}
\newcommand{\props}{\mathsf{Prop}}

\DeclareMathOperator{\free}{F}
\DeclareMathOperator{\forget}{U}
\DeclareMathOperator{\yoneda}{\mathcal{Y}}

\newcommand{\loops}{\Omega}
\DeclareMathOperator{\cobar}{\mathbf{\Omega}}

\DeclareMathOperator{\Sq}{Sq}
\DeclareMathOperator{\ind}{ind}

\DeclareMathOperator{\chains}{N}

\DeclareMathOperator{\gchains}{C}

\DeclarePairedDelimiter\bars{\lvert}{\rvert}

\DeclarePairedDelimiter\angles{\langle}{\rangle}
\DeclarePairedDelimiter\set{\{}{\}}
\DeclarePairedDelimiter\ceil{\lceil}{\rceil}
\DeclarePairedDelimiter\floor{\lfloor}{\rfloor}

\newcommand{\id}{\mathsf{id}}
\renewcommand{\th}{\mathrm{th}}
\newcommand{\op}{\mathrm{op}}

\newcommand{\Hom}{\mathrm{Hom}}
\newcommand{\End}{\mathrm{End}}
\newcommand{\coEnd}{\mathrm{coEnd}}
\newcommand{\xla}[1]{\xleftarrow{#1}}
\newcommand{\xra}[1]{\xrightarrow{#1}}

\newcommand{\pdfEinfty}{\texorpdfstring{${E_\infty}$}{E-infty}}

\newcommand{\rH}{\mathrm{H}}

\newcommand{\cM}{\mathcal{M}}

\newcommand{\cO}{\mathcal{O}}
\newcommand{\cP}{\mathcal{P}}

\newcommand{\cT}{\mathcal{T}}

\newcommand{\sA}{\mathsf{A}}
\newcommand{\sB}{\mathsf{B}}
\newcommand{\sC}{\mathsf{C}}

\newcommand{\fZ}{\mathfrak{Z}}

\newcommand{\symbimod}{\mathsf{biMod}_{\sym}}
\newcommand{\graphs}{\mathfrak{G}}
\newcommand{\biEnd}{\mathrm{biEnd}}

\newcommand{\M}{\cM}
\newcommand{\UM}{{\forget(\M)}}
\newcommand{\UMsh}{{\forget(\M_{sh})}}

\DeclareMathOperator{\schains}{N^{\simplex}}
\DeclareMathOperator{\cchains}{N^{\cube}}
\DeclareMathOperator{\scochains}{N^{\bullet}_{\simplex}}
\DeclareMathOperator{\ccochains}{N^{\bullet}_{\cube}}
\DeclareMathOperator{\Schains}{S}
\DeclareMathOperator{\sSchains}{\Schains^{\simplex}}
\DeclareMathOperator{\cSchains}{\Schains^{\cube}}
\DeclareMathOperator{\sScochains}{S^{\bullet}_{\simplex}}
\DeclareMathOperator{\cScochains}{S^{\bullet}_{\cube}}
\DeclareMathOperator{\sSing}{Sing^{\simplex}}
\DeclareMathOperator{\cSing}{Sing^{\cube}}

\newcommand{\scube}[1]{(\triangle^{\!1})^{\times #1}}
\DeclareMathOperator{\triangulate}{\mathcal{T}}
\DeclareMathOperator{\cubify}{\mathcal{U}}
\DeclareMathOperator{\ez}{\mathfrak{ez}}

\DeclareMathOperator{\cs}{\mathfrak{cs}}
\DeclareMathOperator{\CS}{CS}

\DeclareMathOperator{\ccs}{\mathsf{cs}}

\DeclareMathOperator{\gi}{\mathfrak{i}}

\title[A combinatorial $E_\infty$-algebra structure on cubical cochains]{A combinatorial $E_\infty$-algebra structure on cubical cochains and the Cartan--Serre map}

\author[R.~Kaufmann]{Ralph~M.~Kaufmann}
\address{R.K., Department of Mathematics, Department of Physics and Astronomy, Purdue University}
\email{\href{mailto:rkaufman@purdue.edu}{rkaufman@purdue.edu}}

\author[A.~Medina-Mardones]{Anibal~M.~Medina-Mardones}
\address{A.M-M., Max Plank Institute for Mathematics \and University of Notre Dame}
\email{\href{mailto:ammedmar@mpim-bonn.mpg.de}{ammedmar@mpim-bonn.mpg.de}}

\date{\today}

\begin{document}
	\usetikzlibrary{arrows,decorations.markings}
\tikzset{myptr/.style={decoration={markings,mark=at position 1 with %
			{\arrow[scale=2,>=stealth]{>}}},postaction={decorate}}}
		
\newsavebox\preproduct
\begin{lrbox}{\preproduct}
	\begin{tikzpicture}[scale=.3]
	\draw (0,0)--(0,-.8);
	\draw (0,0)--(.5,.5);
	\draw (0,0)--(-.5,.5);
	\end{tikzpicture} 
\end{lrbox}
\newcommand{\product}{
	\usebox\preproduct}

\newsavebox\precoproduct
	\begin{lrbox}{\precoproduct}
		\begin{tikzpicture}[scale=.3]
		\draw (0,0)--(0,.8);
		\draw (0,0)--(.5,-.5);
		\draw (0,0)--(-.5,-.5);
		\end{tikzpicture}
	\end{lrbox}
\newcommand{\coproduct}{
	\usebox\precoproduct}

\newsavebox\preboundary
\begin{lrbox}{\preboundary}
	\begin{tikzpicture}[scale=.3]
	\draw (0,0)--(0,1.3);
	\draw (.5,0)--(.5,1.3);
	\draw [fill] (0,0) circle [radius=0.1];
	\draw (.9,.6)--(1.3,.6);
	\draw (1.7,0)--(1.7,1.3);
	\draw (2.2,0)--(2.2,1.3);
	\draw [fill] (2.2,0) circle [radius=0.1];
	\end{tikzpicture}
\end{lrbox}
\newcommand{\boundary}{
	\usebox\preboundary}

\newsavebox\precoboundary
\begin{lrbox}{\precoboundary}
	\begin{tikzpicture}[scale=.3]
	\draw (0,0)--(0,1.3);
	\draw (.5,0)--(.5,1.3);
	\draw [fill] (0,1.3) circle [radius=0.1];
	\draw (.9,.6)--(1.3,.6);
	\draw (1.7,0)--(1.7,1.3);
	\draw (2.2,0)--(2.2,1.3);
	\draw [fill] (2.2,1.3) circle [radius=0.1];
	\end{tikzpicture}
\end{lrbox}
\newcommand{\coboundary}{
	\usebox\precoboundary}

\newsavebox\precounit
\begin{lrbox}{\precounit}
	\begin{tikzpicture}[scale=.3]
	\draw (0,0)--(0,1.3);
	\draw [fill] (0,0) circle [radius=0.1];
	\end{tikzpicture}
\end{lrbox}
\newcommand{\counit}{
	\usebox\precounit}

\newsavebox\preidentity
\begin{lrbox}{\preidentity}
	\begin{tikzpicture}[scale=.3]
	\draw (0,0)--(0,1.3);
	\end{tikzpicture}
\end{lrbox}
\newcommand{\identity}{
	\usebox\preidentity}

\newsavebox\preunit
\begin{lrbox}{\preunit}
	\begin{tikzpicture}[scale=.3]
	\draw (0,0)--(0,1.3);
	\draw [fill] (0,1.3) circle [radius=0.1];
	\end{tikzpicture}
\end{lrbox}
\newcommand{\unit}{
	\usebox\preunit}

\newsavebox\preassociativity
\begin{lrbox}{\preassociativity}
	\begin{tikzpicture}[scale=.2]
	\path[draw] (0,0)--(0,1)--(-1,2);
	\draw (0,1)--(1,2);
	\draw (-.5,1.5)--(0,2);
	\draw (1.2,1)--(1.8,1);
	\path[draw] (3,0)--(3,1)--(2,2);
	\draw (3,1)--(4,2);
	\draw (3.5,1.5)--(3,2);	
	\end{tikzpicture}
\end{lrbox}
\newcommand{\associativity}{
	\usebox\preassociativity}

\newsavebox\precoassociativity
\begin{lrbox}{\precoassociativity}
	\begin{tikzpicture}[scale=.2]
	\path[draw] (0,0)--(0,-1)--(-1,-2);
	\draw (0,-1)--(1,-2);
	\draw (-.5,-1.5)--(0,-2);
	\draw (1.2,-1)--(1.8,-1);
	\path[draw] (3,0)--(3,-1)--(2,-2);
	\draw (3,-1)--(4,-2);
	\draw (3.5,-1.5)--(3,-2);
	\end{tikzpicture}
\end{lrbox}
\newcommand{\coassociativity}{
	\usebox\precoassociativity}

\newsavebox\preinvolution
\begin{lrbox}{\preinvolution}
	\begin{tikzpicture}[scale=.2]
	\path[draw] (0,0)--(0,.5)--(-.5,1)--(0,1.5)--(0,2);
	\path[draw] (0,.5)--(.5,1)--(0,1.5);
	\end{tikzpicture}
\end{lrbox}
\newcommand{\involution}{
	\usebox\preinvolution}

\newsavebox\preleftcounitality
\begin{lrbox}{\preleftcounitality}
	\begin{tikzpicture}[scale=.3]
	\draw (0,0)--(0,.8);
	\draw (0,0)--(.5,-.5);
	\draw (0,0)--(-.5,-.5);
	\draw [fill] (-.5,-.5) circle [radius=0.1];
	\draw (.7,0)--(1.1,0);
	\path[draw] (1.5,-.5)--(1.5,.8);
	\end{tikzpicture}
\end{lrbox}
\newcommand{\leftcounitality}{
	\usebox\preleftcounitality}

\newsavebox\preleftcounitcoproduct
\begin{lrbox}{\preleftcounitcoproduct}
	\begin{tikzpicture}[scale=.3]
	\draw (0,0)--(0,.8);
	\draw (0,0)--(.5,-.5);
	\draw (0,0)--(-.5,-.5);
	\draw [fill] (-.5,-.5) circle [radius=0.1];
	\end{tikzpicture}
\end{lrbox}
\newcommand{\leftcounitcoproduct}{
	\usebox\preleftcounitcoproduct}

\newsavebox\prerightcounitality
\begin{lrbox}{\prerightcounitality}
	\begin{tikzpicture}[scale=.3]
	\draw (0,0)--(0,.8);
	\draw (0,0)--(.5,-.5);
	\draw (0,0)--(-.5,-.5);
	\draw [fill] (.5,-.5) circle [radius=0.1];
	\draw (-.7,0)--(-1.1,0);
	\path[draw] (-1.5,-.5)--(-1.5,.8);
	\end{tikzpicture}
\end{lrbox}
\newcommand{\rightcounitality}{
	\usebox\prerightcounitality}

\newsavebox\prerightcounitcoproduct
\begin{lrbox}{\prerightcounitcoproduct}
	\begin{tikzpicture}[scale=.3]
	\draw (0,0)--(0,.8);
	\draw (0,0)--(.5,-.5);
	\draw (0,0)--(-.5,-.5);
	\draw [fill] (.5,-.5) circle [radius=0.1];
	\end{tikzpicture}
\end{lrbox}
\newcommand{\rightcounitcoproduct}{
	\usebox\prerightcounitcoproduct}

\newsavebox\preleftunitality
\begin{lrbox}{\preleftunitality}
	\begin{tikzpicture}[scale=.3]
	\draw (0,0)--(0,-.8);
	\draw (0,0)--(-.5,.5);
	\draw (0,0)--(.5,.5);
	\draw [fill] (.5,.5) circle [radius=0.1];
	\draw (.7,0)--(1.1,0);
	\path[draw] (1.5,.5)--(1.5,-.8);
	\end{tikzpicture}
\end{lrbox}
\newcommand{\leftunitality}{
	\usebox\preleftunitality}

\newsavebox\prerightunitality
\begin{lrbox}{\prerightunitality}
	\begin{tikzpicture}[scale=.3]
	\draw (0,0)--(0,-.8);
	\draw (0,0)--(-.5,.5);
	\draw (0,0)--(.5,.5);
	\draw [fill] (-.5,.5) circle [radius=0.1];
	\draw (-.7,0)--(-1.1,0);
	\path[draw] (-1.5,.5)--(-1.5,-.8);
	\end{tikzpicture}
\end{lrbox}
\newcommand{\rightunitality}{
	\usebox\prerightunitality}

\newsavebox\preproductcounit
\begin{lrbox}{\preproductcounit}
	\begin{tikzpicture}[scale=.3]
	\draw (0,0)--(0,-.8);
	\draw (0,0)--(.5,.5);
	\draw (0,0)--(-.5,.5);
	\draw [fill] (0,-.8) circle [radius=0.1];
	\end{tikzpicture}
\end{lrbox}
\newcommand{\productcounit}{
	\usebox\preproductcounit}

\newsavebox\preunitcoproduct
\begin{lrbox}{\preunitcoproduct}
	\begin{tikzpicture}[scale=.3]
	\draw (0,0)--(0,.8);
	\draw (0,0)--(.5,-.5);
	\draw (0,0)--(-.5,-.5);
	\draw [fill] (0,.8) circle [radius=0.1];
	\end{tikzpicture}
\end{lrbox}
\newcommand{\unitcoproduct}{
	\usebox\preunitcoproduct}

\newsavebox\preleibniz
\begin{lrbox}{\preleibniz}
	\begin{tikzpicture}[scale=.245]
	\draw (0,.3)--(0,-.3);
	\draw (0,.3)--(.5,.8);
	\draw (0,.3)--(-.5,.8);
	\draw (0,-.3)--(0,.3);
	\draw (0,-.3)--(.5,-.8);
	\draw (0,-.3)--(-.5,-.8);
	
	\draw (.7,0)--(1.1,0);
	\draw (2.1,.8)--(2.1,.3)--(1.7,0)--(1.7,-.8);
	\draw (2.1,.3)--(2.9,-.3);
	\draw (3.3,.8)--(3.3,0)--(2.9,-.3)--(2.9,-.8);
	
	\draw (3.8,0)--(4.1,0);
	\draw (4.6,.8)--(4.6,0)--(5,-.3)--(5,-.8);
	\draw (5,-.3)--(5.8,.3);
	\draw (5.8,.8)--(5.8,.3)--(6.2,0)--(6.2,-.8);	
	\end{tikzpicture}
\end{lrbox}
\newcommand{\leibniz}{
	\usebox\preleibniz}

\newsavebox\prebialgebra
\begin{lrbox}{\prebialgebra}
	\begin{tikzpicture}[scale=.5]
	\draw (0,.3)--(0,-.3);
	\draw (0,.3)--(.5,.8);
	\draw (0,.3)--(-.5,.8);
	\draw (0,-.3)--(0,.3);
	\draw (0,-.3)--(.5,-.8);
	\draw (0,-.3)--(-.5,-.8);
	
	\draw (.8,0)--(1.2,0);
	
	\draw (2.1,.8)--(2.1,.3)--(1.7,0)--(2.1,-.3)--(2.1,-.8);
	
	\draw (3.1,.8)--(3.1,.3)--(3.5,0)--(3.1,-.3)--(3.1,-.8);

	\draw (2.1,.3)--(3.1,-.3);
	\draw (2.1,-.3)--(2.5,-.06);
	\draw (3.1,.3)--(2.7,.06);	
	\end{tikzpicture}
\end{lrbox}
\newcommand{\bialgebra}{
	\usebox\prebialgebra}

\newsavebox\precommutativity
\begin{lrbox}{\precommutativity}
	\begin{tikzpicture}[scale=.27]
	\draw (.3,0)--(.3,-1);
	\draw (.3,0)--(.8,.5);
	\draw (.3,0)--(-.2,.5);
	
	\draw (1.2,-.4)--(1.8,-.4);
	
	\draw (3,-.3)--(3,-1);
	\draw (2.5,0)--(3,-.3);
	\draw (3.5,0)--(3,-.3);
	\draw (2.46,0)--(2.9,.18);
	\draw (3.1,.3)--(3.5,.5);
	\draw (3.5,0)--(2.5,.5);
		\end{tikzpicture} 
	\end{lrbox}
	\newcommand{\commutativity}{
		\usebox\precommutativity}	
	\begin{abstract}
	Cubical cochains are equipped with an associative product, dual to the Serre diagonal, lifting the graded commutative structure in cohomology.
	In this work we introduce through explicit combinatorial methods an extension of this product to a full $E_\infty$-structure.
	As an application we prove that the Cartan--Serre map, which relates the cubical and simplicial singular cochains of spaces, is a quasi-isomorphism of $E_\infty$-algebras.
\end{abstract}
	\maketitle
	\section{Introduction} \label{s:introduction}

Instead of simplices, in his groundbreaking work on fibered spaces Serre considered cubes as the basic shapes used to define cohomology, stating that:
\begin{displaycquote}[p.431]{serre1951homologie}
	Il est en effet evident que ces derniers se pretent mieux que les simplexes a l'etude des produits directs, et, a fortiori, des espaces fibres qui en sont la generalisation.
\end{displaycquote}
Cubical sets, a model for the homotopy category, were considered by Kan \cite{kan1955abstract, kan1956abstract} before introducing simplicial sets, they are central to nonabelian algebraic topology \cite{brown2011nonabelian}, and have become important in Voevodsky's program for univalent foundations and homotopy type theory \cite{kapulkin2020straightening, mortberg2017cubical}.
Other areas that highlight the relevance of cubical methods are applied topology, where cubical complexes are ubiquitous in the study of images \cite{tomasz2004computational}, condensed matter physics, where models on cubical lattices are central \cite{baxter1985exactlysolved}, and geometric group theory \cite{gromov1987hyperbolic}, where fundamental results have been obtained considering actions on certain cube complexes characterized combinatorially \cite{agol2013haken}.

Cubical cochains are equipped with the \textit{Serre algebra structure}, a lift to the cochain level of the graded ring structure in cohomology.
Using an acyclic carrier argument it can be shown that this product is commutative up to coherent homotopies in a non-canonical way.
The study of such objects, referred to as $E_\infty$-algebras, has a long history, where (co)homology operations \cite{steenrod1962cohomology, may1970general}, the recognition of infinite loop spaces \cite{boardman1973homotopy, may1972geometry} and complete algebraic models of the $p$-adic homotopy category \cite{mandell2001padic} are key milestones.
The goal of this work is to introduce a description of an explicit $E_\infty$-algebra structure naturally extending the Serre algebra structure, and relate it to one on simplicial cochains extending the Alexander--Whitney algebra structure.

We use the combinatorial model of the $E_\infty$-operad $\UM$ obtained from the finitely presented prop $\M$ introduced in \cite{medina2020prop1}.
The resulting $\UM$-algebra structure on cubical cochains is induced from a natural $\M$-bialgebra structure on the chains of representable cubical sets, which is determined by only three linear maps.
To our knowledge, this is the first effective construction of an $E_\infty$-algebra structure on cubical cochains.
Non-constructively, this result could be obtained using a lifting argument based on the cofibrancy of the reduced version of the operad $\UM$ in the model category of operads \cite{hinich1997homological, berger2003modelcategory}, but this existence statement is not very useful in concrete situations.
To illustrate the advantages of an effective construction let us consider a prime $p$.
The mod $p$ cohomology of spaces is equipped with natural stable endomorphisms, known as Steenrod operations \cite{steenrod1962cohomology}.
Following an operadic viewpoint developed by May \cite{may1970general}, in \cite{medina2021may_st} we exhibited integral elements in $\UM$ representing Steenrod operations on the mod~$p$ homology of $\UM$-algebras.
Since, as proven in this article, the cochains of a cubical set are equipped with a $\UM$-algebra structure, we obtain natural cochain level multioperations for cubical sets representing Steenrod operation at every $p$.
This cubical cup-$(p,i)$ products are explicit enough to have been implemented in the open source computer algebra system \href{https://comch.readthedocs.io/en/latest/}{\texttt{ComCH}} \cite{medina2021comch}.

We now turn to the comparison between cubical and simplicial cochains.
In \cite[p. 442]{serre1951homologie}, Serre described for any topological space $\fZ$ a natural quasi-isomorphism
\begin{equation} \label{e:cs on singular cochains}
	\cScochains(\fZ) \to \sScochains(\fZ)
\end{equation}
between its cubical and simplicial singular cochains, stating this to be a quasi-isomorphism of algebras with respect to the usual structures.
We will consider a well known Quillen equivalence
\[
\begin{tikzcd}[column sep=0]
	\sSet \arrow[rr, "\cubify"', bend right] & \perp & \arrow[ll,"\triangulate"', bend right] \cSet
\end{tikzcd}
\]
between simplicial and cubical sets, and construct a natural chain map
\begin{equation} \label{e:intro main map}
	\ccochains(\cubify Y) \to \scochains(Y)
\end{equation}
for every simplicial set $Y$.
In \cite{medina2020prop1}, a natural $\UM$-algebra structure extending the Alexander--Whitney coalgebra structure was constructed on simplicial sets.
With respect to it and the one defined here for cubical sets we have the following results after passing to a sub-$E_\infty$-operad of $\UM$.

\begin{theorem*}
	The map presented in \cref{e:intro main map} is a quasi-isomorphism of $E_\infty$-algebras.
\end{theorem*}

From this result, stated as \cref{t:main comparison}, we deduce the following two.
The first one concerns the triangulation functor $\triangulate$ and it is stated more precisely as \cref{c:zig-zag}.

\begin{corollary*}
	There is a natural zig-zag of $E_\infty$-algebra quasi-isomorphisms between the cochains of a cubical set and those of its triangulation.
\end{corollary*}

The next one concerns the map presented in \cref{e:cs on singular cochains}, relating the cubical and simplicial singular cochains of a space, and it is stated more precisely as \cref{c:cs e infty}.

\begin{corollary*}
	The Cartan--Serre map is a quasi-isomorphism of $E_\infty$-algebras.
\end{corollary*}

\begin{remark*}
	In this introduction we have used the setting defined by cochains and products since it is more familiar, whereas in the rest of the text we use the more fundamental one defined by chains and coproducts.
\end{remark*}

\section*{Outline}

We recall the required notions from homological algebra and category theory in \cref{s:preliminaries}.
The necessary concepts from the theory of operads and props is reviewed in \cref{s:props}, including the definition of the prop $\M$.
\cref{s:action} contains our main contribution; an explicit natural $\M$-bialgebra structure on the chains of representable cubical sets and, from it, a natural $E_\infty$-coalgebra structure on the chains of cubical sets.
The comparison between simplicial and cubical chains is presented in \cref{s:comparison}, where we show that the Cartan--Serre map is a quasi-isomorphism respecting $E_\infty$-structures.
We close presenting some future work in \cref{s:future}.
	\section*{Acknowledgment}

We thank the reviewer for many valuable suggestions improving the presentation of this work.
We are grateful to Clemens Berger, Greg Friedman, Kathryn Hess, Peter May, Manuel Rivera, Paolo Salvatore, Dev Sinha, Dennis Sullivan, Peter Teichner, and Bruno Vallette for insightful discussions related to this project.
A.M-M. acknowledges financial support from Innosuisse grant 32875.1~IP-ICT-1 and the hospitality of the \textit{Laboratory for Topology and Neuroscience} at EPFL.
	\section{Conventions and preliminaries} \label{s:preliminaries}

\subsection{Chain complexes}

Throughout this article $\k$ denotes a commutative and unital ring and we work over its associated closed symmetric monoidal category of differential (homologically) graded $\k$-modules $(\Ch, \ot, \k)$.
We refer to the objects and morphisms of this category as \textit{chain complexes} and \textit{chain maps} respectively.
We denote by $\Hom(C, C^\prime)$ the chain complex of $\k$-linear maps between chain complexes $C$ and $C^\prime$, and refer to the functor $\Hom(-, \k)$ as \textit{linear duality}.

\subsection{Presheaves}

Recall that a category is said to be \textit{small} if its objects and morphisms form sets.
We denote the category of small categories by $\Cat$.
Given categories $\sB$ and $\sC$ with $\sB$ small we denote their associated \textit{functor category} by $\Fun(\sB, \sC)$.
A category is said to be \textit{cocomplete} if any functor to it from a small category has a colimit.
If $\sA$ is small and $\sC$ cocomplete, then the (\textit{left}) \textit{Kan extension of $g$ along $f$} exists for any pair of functors $f$ and $g$ in the diagram below, and it is the initial object in $\Fun(\sB, \sC)$ making
\[
\begin{tikzcd}[column sep=normal, row sep=normal]
\sA \arrow[d, "f"'] \arrow[r, "g"] & \sC \\
\sB \arrow[dashed, ur, bend right] & \quad
\end{tikzcd}
\]
commute.
A Kan extension along the \textit{Yoneda embedding}, i.e., the functor
\[
\yoneda \colon \sA \to \Fun(\sA^\op, \Set)
\]
induced by the assignment
\[
a \mapsto \big( a^\prime \mapsto \sA(a^\prime, a) \big),
\]
is referred to as a \textit{Yoneda extension}.
Abusively we use the same notation for a functor and for its Yoneda extension.
We refer to objects of $\Fun(\sA^\op, \Set)$ in the image of the Yoneda embedding as \textit{representable}.
	\section{Operads, props and \pdfEinfty-structures} \label{s:props}

We now review the definition of the finitely presented prop $\M$ introduced in \cite{medina2020prop1} and whose associated operad is a model of the $E_\infty$-operad.
Given its small number of generators and relations, is well suited to explicitly define $E_\infty$-structures.
We start by recalling some basic material from the theory of operads and props.

\subsection{Symmetric (bi)modules}

Let $\sym$ be the category whose objects are the non-negative integers $\N$ and whose set of morphisms between $n$ and $n'$ is empty if $n \neq n'$ and is otherwise the symmetric group $\sym_n$.
A \textit{left $\sym$-module} (resp. \textit{right} $\sym$-\textit{module} or $\sym$-\textit{bimodule}) is a functor from $\sym$ (resp. $\sym^\op$ or $\sym^\op \times \sym$) to $\Ch$.
In this paper we prioritize left module structures over their right counterparts.
As usual, taking inverses makes both perspectives equivalent.
We respectively denote by $\symMod$ and $\symbimod$ the categories of left $\sym$-modules and of $\sym$-bimodules with morphisms given by natural transformations.

Given a chain complex $C$, we have the following key examples of a left and a right $\sym$-module
\begin{gather*}
	\End^C(n) = \Hom(C, C^{\ot n}), \qquad
	\End_C(m) = \Hom(C^{\ot m}, C),
\end{gather*}
and of an $\sym$-bimodule
\[
\End^C_C(m,n) = \Hom(C^{\ot m}, C^{\ot n}),
\]
where the symmetric actions are given by permutation of tensor factors.

The group homomorphisms $\sym_n \to \sym_1^\op \times \sym_n$ induce a forgetful functor
\begin{equation} \label{e:forgetful}
	\forget \colon \symbimod \to \symMod
\end{equation}
defined explicitly on an object $\cP$ by $\forget(\cP)(n) = \cP(1, n)$ for $n \in \N$.
The similarly defined forgetful functor to right $\sym$-modules will not be considered.

\subsection{Composition structures}

\textit{Operads} and \textit{props} are obtained by enriching $\sym$-modules and $\sym$-bimodules with certain composition structures.
Intuitively, these are obtained by abstracting the composition structure naturally present in the left $\sym$-module $\End^C$ (or right $\sym$-module $\End_C$), naturally an operad, and the $\sym$-bimodule $\End^C_C$, naturally a prop.
More explicitly, an operad $\cO$ is a left $\sym$-module with chain maps
\begin{gather*}
	\k \to \cO(1), \\
	\cO(n_1) \ot \dotsb \ot \cO(n_r) \ot \cO(r) \to \cO(n_1 + \dots + n_r),
\end{gather*}
satisfying relations of associativity, equivariance and unitality.
Similarly, a prop $\cP$ is an $\sym$-bimodule together with chain maps
\begin{gather*}
	\k \to \cP(n,n), \\
	\cP(m,k) \ot \cP(k,n) \to \cP(m,n), \\
	\cP(m,n) \ot \cP(m',n') \to \cP(m+m',n+n'),
\end{gather*}
satisfying certain natural relations.
For a complete presentation of these concepts we refer to Definition~11 and 54 of \cite{markl2008props}.
We respectively denote the category of operads and props with structure preserving morphisms by $\operads$ and $\props$.

Let $C$ be a chain complex, $\cO$ an operad, and $\cP$ a prop.
An $\cO$-\textit{coalgebra} (resp. $\cO$-\textit{algebra} or $\cP$-\textit{bialgebra}) structure on $C$ is a structure preserving morphism $\cO \to \End^C$ (resp. $\cO \to \End_C$ or $\cP \to \End_C^C$).
We mention that the linear dual of an $\cO$-coalgebra is an $\cO$-algebra.

Since the forgetful functor presented in \cref{e:forgetful} induces a functor
\[
\forget \colon \props \to \operads,
\]
any $\cP$-bialgebra structure on $C$
\[
\cP \to \biEnd_C^C
\]
induces a $\forget(\cP)$-coalgebra structure on it
\[
\forget(\cP) \to \forget(\biEnd_C^C) \cong \coEnd^C.
\]

\subsection{\pdfEinfty-operads}

Recall that a \textit{projective $\sym_n$-resolution} of a chain complex $C$ is a quasi-isomorphism $R \to C$ from a chain complex $R$ of projective $\k[\sym_n]$-modules.
An $\sym$-module $M$ is said to be $E_\infty$ if there exists a morphism of $\sym$-modules $M \to \underline{\k}$ inducing for each $n \in \N$ a free $\sym_n$-resolution $M(n) \to \k$.
An operad is said to be an \textit{$E_\infty$-operad} if its underlying $\sym$-module is $E_\infty$.
A prop $\cP$ is said to be an \textit{$E_\infty$-prop} if $\forget(\cP)$ is an $E_\infty$-operad.

\subsection{Presentations} \label{ss:presentation}

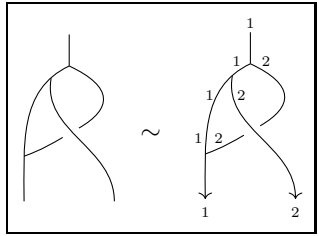
\begin{figure}[h]
	\centering
	\boxed{
	\begin{tikzpicture}[scale=.6]
		\draw (1,3.7) to (1,3);

		\draw (1,3) to [out=205, in=90] (0,0);

		\draw [shorten >= 0cm] (.6,2.73) to [out=-100, in=90] (2,0);

		\draw [shorten >= .15cm] (1,3) to [out=-25, in=30, distance=1.1cm] (1,1.5);
		\draw [shorten <= .1cm] (1,1.5) to [out=210, in=20] (0,1);

		\node at (1,3.9){};
		\node at (0,-.32){};
		\node at (2,-.32){};

		\node at (3,1.5){$\sim$\ \ \ };
	\end{tikzpicture}
	\begin{tikzpicture}[scale=.6]
		\draw (1,3.7) to (1,3);

		\draw [->](1,3) to [out=205, in=90] (0,0);

		\draw [shorten >= 0cm,->] (.6,2.73) to [out=-100, in=90] (2,0);

		\draw [shorten >= .15cm] (1,3) to [out=-25, in=30, distance=1.1cm] (1,1.5);
		\draw [shorten <= .1cm] (1,1.5) to [out=210, in=20] (0,1);

		\def\x{.8}

		\node[scale=\x] at (1,3.9){$\scriptstyle 1$};

		\node[scale=\x] at (.7,3.05){$\scriptstyle 1$};
		\node[scale=\x] at (1.35,3.05){$\scriptstyle 2$};

		\node[scale=\x] at (.1,2.3){$\scriptstyle 1$};
		\node[scale=\x] at (.8,2.3){$\scriptstyle 2$};

		\node[scale=\x] at (-.15,1.35){$\scriptstyle 1$};
		\node[scale=\x] at (.3,1.35){$\scriptstyle 2$};

		\node[scale=\x] at (0,-.3){$\scriptstyle 1$};
		\node[scale=\x] at (2,-.3){$\scriptstyle 2$};
	\end{tikzpicture}
}
	\caption{Immersed graph representing a $(1,2)$-graph.}
	\label{f:immersion}
\end{figure}

The \textit{free prop} construction is the left adjoint to the forgetful functor from props to $\sym$-bimodules.
Explicitly, the free prop $\free(M)$ generated by an $\sym$-bimodule $M$ is constructed using isomorphism classes of directed graphs with no directed loops that are enriched with the following labeling structure.
We think of each directed edge as built from two compatibly directed half-edges.
For each vertex $v$ of a directed graph $\Gamma$, we have the sets $in(v)$ and $out(v)$ of half-edges that are respectively incoming to and outgoing from $v$.
Half-edges that do not belong to $in(v)$ or $out(v)$ for any $v$ are divided into the disjoint sets $in(\Gamma)$ and $out(\Gamma)$ of incoming and outgoing external half-edges.
For any positive integer $n$ let $\overline{n} = \{1, \dots, n\}$ and set $\overline{0} = \emptyset$.
For any finite set $S$, denote the cardinality of $S$ by $|S|$.
The labeling is given by bijections
\[
\overline{|in(\Gamma)|} \to in(\Gamma), \qquad
\overline{|out(\Gamma)|} \to out(\Gamma),
\]
and
\[
\overline{|in(v)|} \to in(v), \qquad
\overline{|out(v)|} \to out(v),
\]
for every vertex $v$.
We refer to the isomorphism classes of such labeled directed graphs with no directed loops and $m$ incoming and $n$ outgoing half-edges as $(m,n)$\textit{-graphs}.
We denote the set these form by $\graphs(m,n)$.
We use graphs immersed in the plane to represent elements in $\graphs(m,n)$, with the direction implicitly given from top to bottom and the labeling from left to right.
Please consult \cref{f:immersion} for an example.

We consider the right action of $\sym_m$ and the left action of $\sym_n$ on a $(m,n)$-graph given respectively by permuting the labels of $in(\Gamma)$ and $out(\Gamma)$.
This action defines the $\sym$-bimodule structure on the free prop
\begin{equation} \label{e:free prop}
\free(M)(m,n) \ =
\bigoplus_{\substack{\Gamma \text{ in } \\ \graphs(m,n)}} \
\bigotimes_{\substack{v \text{ in} \\ Vert(\Gamma)}} out(v) \ot_{\sym_p} M(p, q) \ot_{\sym_p} in(v),
\end{equation}
where we have simplified the notation writing $p$ and $q$ for $\overline{|in(v)|}$ and $\overline{|out(v)|}$ respectively.
The differential $\bd_{\free(M)}$ is the extension of that of $M$ to the tensor product \eqref{e:free prop}, and the prop structure is induced by
the ``identity graphs''
\[
\identity \ \identity \dotsb \identity
\]
together with (relabeled) grafting and disjoint union.

Let $G$ be an assignment of a set $G(m,n)_d$ to each $m,n \in \N$ and $d \in \Z$.
Denote by $\k[\sym^\op \times \sym] \set{G}$ the $\sym$-bimodule mapping $(m,n)$ to the chain complex with trivial differential and degree $d$ part equal to
\[
\k[\sym^\op_m \times \sym_n] \set[\big]{G(m,n)_d}.
\]
We will denote by $\free(G)$ the free prop generated by this $\sym$-bimodule.
Let $\bd \colon \k[\sym^\op \times \sym] \set{G} \to \free(G)$ be a morphism of $\sym$-bimodules whose canonical extension $\bd \colon \free(G) \to \free(G)$ defines a differential.
We denote by $\free_{\bd}(G)$ the prop obtained by endowing $\free(G)$ with this differential.
Let $R$ be a collection of elements in $\free(G)$ and denote by $\angles{R\hspace{1.5pt}}$ the smallest ideal containing $R$.
The prop \textit{generated by $G$ modulo $R$ with boundary $\bd$} is defined to be $\free_{\bd}(G) \big/ \angles{R\hspace{1.5pt}}$.

\subsection{The prop $\M$}

We now recall the $E_\infty$-prop that is central to our constructions.

\begin{definition}
	Let $\M$ be the prop generated by
	\begin{equation} \label{e:generators of M}
	\counit \,, \qquad
	\coproduct \,, \qquad
	\product,
	\end{equation}
	in $(1,0)_0$\,, $(1,2)_0$ and $(2,1)_1$ respectively,
	modulo the relations
	\begin{equation} \label{e:relations of M}
		\leftcounitality \,, \qquad
		\rightcounitality \,, \qquad
		\productcounit,
	\end{equation}
	with boundary defined by
	\begin{equation} \label{e:boundary of M}
	\bd\ \counit = 0, \qquad
	\bd\, \coproduct = 0, \qquad
	\bd \product = \ \boundary\,.
	\end{equation}
\end{definition}

Explicitly, any element in $\M(m,n)$ can be written as a linear combination of the $(m,n)$-graphs generated by those in \eqref{e:generators of M} via grafting, disjoint union and relabeling, modulo the ideal generated by the relations in \eqref{e:relations of M}. Its boundary is determined by \eqref{e:boundary of M} using \eqref{e:free prop}.

\begin{proposition}[{\cite[Theorem 3.3]{medina2020prop1}}]
	$\M$ is an $E_\infty$-prop.
\end{proposition}

\begin{remark*}
	The prop $\M$ is obtained from applying the functor of chains to a prop over the category of cellular spaces \cite{medina2021prop2}, a quotient of which is isomorphic to the $E_\infty$-operad of stable arc surfaces \cite{kaufmann2009dimension}.
\end{remark*}

	\section{An \pdfEinfty-structure on cubical chains} \label{s:action}

In this section we construct a natural $\M$-bialgebra structure on the chains of representable cubical sets.
These are determined by three natural linear maps satisfying the relations defining $\mathcal M$.
A Yoneda extension then provides the chains of any cubical set with a natural $\UM$-coalgebra structure.
We begin by recalling the basics of cubical topology.

\subsection{Cubical sets}

The objects of the \textit{cube category} $\cube$ are the sets $2^n = \{0, 1\}^n$ with $2^0 = \{0\}$ for $n \in \N$, and its morphisms are generated by the \textit{coface} and \textit{codegeneracy} maps
\begin{align*}
\delta_i^\varepsilon & = \mathrm{id}_{2^{i-1}} \times \delta^\varepsilon \times \mathrm{id}_{2^{n-1-i}} \colon 2^{n-1} \to 2^n, \\
\sigma_i & = \mathrm{id}_{2^{i-1}} \times \, \sigma \times \mathrm{id}_{2^{n-i}} \quad \colon 2^{n} \to 2^{n-1},
\end{align*}
where $\varepsilon \in \{0,1\}$ and the functors
\[
\begin{tikzcd} [column sep=16pt]
2^0 \arrow[r, bend left, "\delta^0"] \arrow[r, bend right, "\delta^1"'] & 2^1 \arrow[r, "\sigma"] & 2^0
\end{tikzcd}
\]
are defined by
\[
\delta^0(0) = 0, \qquad \delta^1(0) = 1, \qquad \sigma(0) = \sigma(1) = 0.
\]
More globally, the category $\cube$ is the free strict monoidal category with an assigned internal bipointed object.
We refer to \cite{grandis2003cubical} for a more leisurely exposition and variants of this definition.

We denote by $\mathrm{Dgn}(2^m, 2^n)$ the subset of morphism in $\cube(2^m, 2^n)$ of the form $\sigma_i \circ \tau$ with $\tau \in \cube(2^m, 2^{n+1})$.

The category of \textit{cubical sets} $\Fun(\cube^\op, \Set)$ is denoted by $\cSet$ and the representable cubical set $\yoneda(2^n)$ by $\cube^n$.
For any cubical set $X$ we write, as usual, $X_n$ instead of $X(2^n)$.

\subsection{Cubical topology}

Consider the topological $n$-cube
\[
\gcube^{n} = \big\{ (x_1, \dots, x_n) \mid x_i \in [0,1] \big\}.
\]
The assignment $2^n \to \gcube^n$ defines a functor $\cube \to \Top$ with
\begin{align*}
\delta^\varepsilon_i(x_1, \dots, x_n) &= (x_1, \dots, x_i, \varepsilon, x_{i+1}, \dots x_n), \\
\sigma_i(x_1,\dots,x_n) &= (x_1, \dots, \widehat{x}_i, \dots, x_n).
\end{align*}
Its Yoneda extension is known as \textit{geometric realization}.
It has a right adjoint $\cSing \colon \Top \to \cSet$ referred to as the \textit{cubical singular complex} satisfying
\[
\cSing(\fZ)_n = \Top(\gcube^n, \fZ)
\]
for any topological space $\fZ$.

\subsection{Cubical chains}

The functor of (\textit{normalized}) \textit{chains} $\chains \colon \cSet \to \Ch$ is the Yoneda extension of the functor $\cube \to \Ch$ defined next.
It assigns to an object $2^n$ the chain complex having in degree $m$ the module
\[
\frac{\k\{\cube(2^m, 2^n)\}}{\k\{\mathrm{Dgn}(2^m, 2^n)\}}
\]
and differential induced by
\[
\bd (\id_{2^n}) = \sum_{i=1}^{n} \ (-1)^i \
\big(\delta_i^1 - \delta_i^0 \big).
\]
To a morphism $\tau \colon 2^n \to 2^{n'}$ it assigns the chain map
\[
\begin{tikzcd}[row sep=-3pt, column sep=normal,
/tikz/column 1/.append style={anchor=base east},
/tikz/column 2/.append style={anchor=base west}]
\chains(\cube^n)_m \arrow[r] & \chains(\cube^{n'})_m \\
\big( 2^m \to 2^n \big) \arrow[r, mapsto] & \big( 2^m \to 2^n \stackrel{\tau}{\to} 2^{n'} \big).
\end{tikzcd}
\]
The chain complex $\chains(\cube^n)$ is isomorphic to both: $\chains(\cube^1)^{\ot n}$ and the cellular chains on the topological $n$-cube with its standard CW structure $\gchains(\gcube^n)$.
We use the isomorphism $\chains(\cube^n) \cong \gchains(\gcube^1)^{\ot n}$ when denoting the elements in the basis of $\chains(\cube^n)$ by $x_1 \ot \dotsb \ot x_n$ with $x_i \in \{[0], [0,1], [1]\}$.

For a topological space $\fZ$, the chain complex $\chains(\cSing \fZ)$ is referred to as the \textit{cubical singular chains} of $\fZ$.

\subsection{Serre coalgebra} \label{ss:serre coalgebra}

We now recall the \textit{Serre coalgebra structure}, a natural (counital and coassociative) coalgebra structure on cubical chains.

By a Yoneda extension, to define this structure it suffices to describe it on the chains of representable cubical sets $\chains(\cube^n)$.
For $\chains(\cube^1)$ we have
\begin{align*}
	\epsilon \big([0]\big) &= 1,
	& \Delta \big([0]\big) &= [0] \ot [0], \\
	\epsilon \big([1]\big) &= 1,
	& \Delta \big([1]\big) &= [1] \ot [1], \\
	\epsilon \big([0,1]\big) &= 0,
	& \Delta \big([0,1]\big) &= [0] \ot [0,1] + [0,1] \ot [1].
\end{align*}
The Serre coalgebra structure on a general $\chains(\cube^n)$ is define using the isomorphism $\chains(\cube^n) \cong \chains(\cube^1)^{\ot n}$ and the monoidal structure on the category of coalgebras.
Explicitly, the structure maps are given by the compositions
\[
\epsilon \colon \chains(\cube^1)^{\ot n} \xra{\epsilon^{\ot n}} \k^{\ot n} \to \k
\]
and
\[
\Delta \colon \chains(\cube^1)^{\ot n} \xra{\Delta^{\!\ot n}} \left( \chains(\cube^1)^{\ot 2} \right)^{\ot n} \xra{\sigma_{2n}^{-1}} \left( \chains(\cube^1)^{\ot n} \right)^{\ot 2},
\]
where $\sigma_{2n}$ in $\sym_{2n}$ is the $(n,n)$-shuffle mapping the first and second ``decks'' to odd and even values respectively.
An explicit description of $\sigma_{2n}$ is presented in \cref{e:shuffle perm}.

\begin{remark*}
	Similarly to how the Alexander--Whitney coalgebra can be interpreted geometrically as the sum of all complementary pairs of front and back faces of a simplex, this coproduct is, up to signs, also given by the sum of complementary pairs of front and back faces of a cube.
\end{remark*}

For later reference we record a useful description of the value of $\Delta$ on the top dimensional basis element of $\chains(\cube^n)$.

\begin{lemma} \label{l:coproduct description}
	For any $n \in \N$,
	\[
	\Delta \big( [0,1]^{\ot n} \big) =
	\sum_{\lambda \in \Lambda} (-1)^{\ind \lambda} \
	\Big(x_1^{(\lambda)} \ot \cdots \ot x_n^{(\lambda)}\Big) \ot
	\Big(y_1^{(\lambda)} \ot \cdots \ot y_n^{(\lambda)}\Big),
	\]
	where each $\lambda$ in $\Lambda$ is a map $\lambda \colon \{1,\dots,n\} \to \set{0, 1}$ with $\lambda(i)$ interpreted as
	\begin{align*}
		0:\ \, x_i^{(\lambda)} &= [0,1], &	1: \ \, x_i^{(\lambda)} &= [0], \\
		y_i^{(\lambda)} &= [1],  & y_i^{(\lambda)} &= [0,1],
	\end{align*}
	and $\ind \lambda$ is the cardinality of $\set{i<j \mid \lambda(i) > \lambda(j)}$.
\end{lemma}

\subsection{Degree 1 product}

Let $n \in \N$.
For $x = x_1 \ot \dotsb \ot x_n$ a basis element of $\chains(\cube^n)$ and $\ell \in \set{1,\dots,n}$ we write
\begin{align*}
	x_{<\ell} & = x_1 \ot \dotsb \ot x_{\ell-1}, \\
	x_{>\ell} & = x_{\ell+1} \ot \dotsb \ot x_n,
\end{align*}
with the convention
\[
x_{<1} = x_{>n} = 1 \in \Z.
\]
We define the \textit{product} $\ast \colon \chains(\square^n)^{\ot 2} \to \chains(\square^n)$ by
\begin{multline*}
(x_1 \ot \dotsb \ot x_n) \ast (y_1 \ot \dotsb \ot y_n)
=
(-1)^{|x|} \sum_{i=1}^n x_{<i}\, \epsilon(y_{<i}) \ot x_i \ast y_i \ot \epsilon(x_{>i}) \, y_{>i},
\end{multline*}
where the only non-zero values of $x_i \ast y_i$ are
\[
[0] \ast [1] = [0, 1], \qquad [1] \ast [0] = -[0, 1].
\]

\begin{figure}[h!]
	\centering
		\begin{tikzpicture}[every edge quotes/.append style={auto, text=black}, scale=.4]
		\pgfmathsetmacro{\cubex}{4}
		\pgfmathsetmacro{\cubey}{4}
		\pgfmathsetmacro{\cubez}{4}
		\draw [draw=black, every edge/.append style={draw=black, densely dashed, opacity=.5}, fill=white]
		(0,0,0) coordinate (o) -- ++(-\cubex,0,0) coordinate (a) -- ++(0,-\cubey,0) coordinate (b) edge coordinate [pos=1] (g) ++(0,0,-\cubez) -- ++(\cubex,0,0) coordinate (c) -- cycle
		(o) -- ++(0,0,-\cubez) coordinate (d) -- ++(0,-\cubey,0) coordinate (e) edge (g) -- (c) -- cycle
		(o) -- (a) -- ++(0,0,-\cubez) coordinate (f) edge (g) -- (d) -- cycle;
		\path [every edge/.append style={draw=black, |-|}];

		\draw [blue,fill] (b) circle [radius=3pt];
		\draw[] node at (2.5,-1){$\ast$};
	\end{tikzpicture}
	\begin{tikzpicture}[every edge quotes/.append style={auto, text=black}, scale=.4]
		\pgfmathsetmacro{\cubex}{4}
		\pgfmathsetmacro{\cubey}{4}
		\pgfmathsetmacro{\cubez}{4}
		\draw [draw=black, every edge/.append style={draw=black, densely dashed, opacity=.5}, fill=white]
		(0,0,0) coordinate (o) -- ++(-\cubex,0,0) coordinate (a) -- ++(0,-\cubey,0) coordinate (b) edge coordinate [pos=1] (g) ++(0,0,-\cubez) -- ++(\cubex,0,0) coordinate (c) -- cycle
		(o) -- ++(0,0,-\cubez) coordinate (d) -- ++(0,-\cubey,0) coordinate (e) edge (g) -- (c) -- cycle
		(o) -- (a) -- ++(0,0,-\cubez) coordinate (f) edge (g) -- (d) -- cycle;
		\path [every edge/.append style={draw=black, |-|}];

		\draw [blue, fill] (d) circle [radius=3pt];
		\draw[] node at (2.5,-1){$=$};
	\end{tikzpicture}
	\begin{tikzpicture}[every edge quotes/.append style={auto, text=black}, scale=.4]
		\pgfmathsetmacro{\cubex}{4}
		\pgfmathsetmacro{\cubey}{4}
		\pgfmathsetmacro{\cubez}{4}
		\draw [draw=black, every edge/.append style={draw=black, densely dashed, opacity=.5}, fill=white]
		(0,0,0) coordinate (o) -- ++(-\cubex,0,0) coordinate (a) -- ++(0,-\cubey,0) coordinate (b) edge coordinate [pos=1] (g) ++(0,0,-\cubez) -- ++(\cubex,0,0) coordinate (c) -- cycle
		(o) -- ++(0,0,-\cubez) coordinate (d) -- ++(0,-\cubey,0) coordinate (e) edge (g) -- (c) -- cycle
		(o) -- (a) -- ++(0,0,-\cubez) coordinate (f) edge (g) -- (d) -- cycle;
		\path [every edge/.append style={draw=black, |-|}];

		\draw[blue, thick] (b)--(a)--(f)--(d);
	\end{tikzpicture}
	\caption{Geometric representation of $\big([0]\ot[0]\ot[0] \ast [1]\ot[1]\ot[1]\big)$ where we are using the width-depth-height order.}
	\label{f:product}
\end{figure}
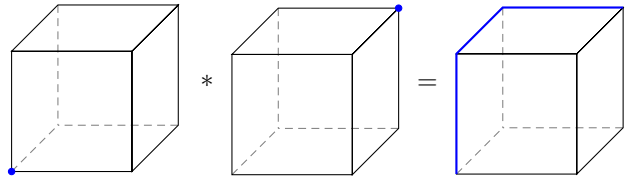

\begin{example*}
	Since in $\chains(\cube^3)$ we have that
	\[
	\bd \big([0]\ot[0]\ot[0]\big) = \bd \big([1]\ot[1]\ot[1]\big) = 0
	\]
	and
	\begin{multline*}
		\bd\big([0]\ot[0]\ot[0] \ast [1]\ot[1]\ot[1]\big) \\ =
		\bd\big([0,1]\ot[1]\ot[1] + [0]\ot[0,1]\ot[1] + [0]\ot[0]\ot[0,1] \big) \\ =
		[1]\ot[1]\ot[1] - [0]\ot[0]\ot[0],
	\end{multline*}
	we conclude that in general $\ast$ is not a cycle in the appropriate Hom complex, so it does not descend to homology.
	This product should be understood as an algebraic version of a consistent choice of path between points in a cube.
	In our case, as illustrated in \cref{f:product}, the chosen path is given by the union of segments parallel to edges of the cube.
\end{example*}

\subsection{$\M$-bialgebra on representable cubical sets}

\begin{lemma} \label{l:cubical chain bialgebra}
	The assignment
	\[
	\counit \mapsto \epsilon, \quad \coproduct \mapsto \Delta, \quad \product \mapsto \ast,
	\]
	induces a natural $\mathcal M$-bialgebra structure on $\chains(\square^n)$ for every $n \in \N$.
\end{lemma}

\begin{proof}
	We need to show that this assignment is compatible with the relations
	\[
	\productcounit = 0, \qquad
	\leftcounitality = 0, \qquad
	\rightcounitality = 0,
	\]
	and
	\[
	\bd\ \counit = 0, \qquad
	\bd\ \coproduct = 0, \qquad
	\bd\ \product = \ \boundary\,.
	\]
	For the rest of this proof let us consider two basis elements of $\chains(\square^n)$
	\begin{align*}
		x = x_1 \ot \dotsb \ot x_n
		\qquad \text{ and } \qquad
		y = y_1 \ot \dotsb \ot y_n.
	\end{align*}
	Since the degree of $\ast$ is $1$ and $\epsilon([0,1]) = 0$, we can verify the first relation easily:
	\begin{align*}
		\varepsilon(x \ast y) & =
		\sum (-1)^{|x|} \epsilon(y_{<i}) \epsilon(x_{<i}) \ot \epsilon(x_i \ast y_i) \ot \epsilon(x_{>i}) \epsilon(y_{>i}) = 0.
	\end{align*}
	For the second relation we want to show that $(\epsilon \ot \id) \circ \Delta = \id$.
	Since
	\begin{gather*}
		(\epsilon \ot \id) \circ \Delta([0]) = \epsilon([0]) \ot [0] = [0], \\
		(\epsilon \ot \id) \circ \Delta([1]) = \epsilon([1]) \ot [1] = [1], \\
		(\epsilon \ot \id) \circ \Delta([0, 1]) = \epsilon([0]) \ot [0, 1] + \epsilon([0, 1]) \ot [1] = [0,1],
	\end{gather*}
	we have
	\begin{multline*}
		(\epsilon \ot \id) \circ \Delta (x_1 \ot \dotsb \ot x_n) = \\
		\sum \pm \left( \epsilon \big(x_1^{(1)}\big) \ot \dotsb \ot \epsilon\big(x_n^{(1)}\big) \right) \ot
		\left( x_1^{(2)} \ot \dotsb \ot x_n^{(2)} \right) \\ =
		x_1 \ot \dotsb \ot x_n,
	\end{multline*}
	where the sign is obtained by noticing that the only non-zero term occurs when each factor $x_i^{(0)}$ is of degree $0$.
	The third relation is verified analogously.
	The fourth and fifth are precisely the well known facts that $\epsilon$ and $\Delta$ are chain maps.
	To verify the sixth and final relation we need to show that
	\[
	\bd (x \ast y)\ +\ \bd x \ast y\ +\ (-1)^{|x|}x \ast \bd y\ =\ \epsilon(x) y \ -\ \epsilon(y) x.
	\]
	We have
	\[
	x \ast y = \sum (-1)^{|x|} x_{<i} \, \epsilon(y_{<i}) \ot x_i \ast y_i \ot \epsilon(x_{>i})\, y_{>i}
	\]
	and
	\begin{align*}
		\bd(x \ast y) & =
		\sum (-1)^{|x|} \, \bd x_{<i}\, \epsilon(y_{<i}) \ot x_i \ast y_i \ot \epsilon(x_{>i})\, y_{>i} \\ & +
		\sum (-1)^{|x|+|x_{<i}|} \, x_{<i}\, \epsilon(y_{<i}) \ot \bd (x_i \ast y_i) \ot \epsilon(x_{>i}) \, y_{>i} \\ & -
		\sum (-1)^{|x|+|x_{<i}|} \, x_{<i}\, \epsilon(y_{<i}) \ot x_i \ast y_i \ot \epsilon(x_{>i})\, \bd y_{>i}.
	\end{align*}
	Since $|x| = |x_{<i}| + |x_i| + |x_{>i}|$ and $\epsilon(x_{>i}) \neq 0 \Leftrightarrow |x_{>i}| = 0$ as well as $\bd(x_i \ast y_i) \neq 0 \Rightarrow |x_i| = 0$
	we have
	\begin{equation} \label{e:boundary of product 1}
		\begin{split}
			\bd(x \ast y) & =
			\sum (-1)^{|x|} \, \bd x_{<i}\, \epsilon(y_{<i}) \ot x_i \ast y_i \ot \epsilon(x_{>i})\, y_{>i} \\ & +
			\sum x_{<i} \, \epsilon(y_{<i}) \ot \bd (x_i \ast y_i) \ot \epsilon(x_{>i})\, y_{>i} \\ & -
			\sum x_{<i} \, \epsilon(y_{<i}) \ot x_i \ast y_i \ot \epsilon(x_{>i})\, \bd y_{>i}.
		\end{split}
	\end{equation}
	We also have
	\begin{align*}
		\bd x \ast y & =
		\sum (-1)^{|x|-1} \, \bd x_{<i}\, \epsilon(y_{<i}) \ot x_i \ast y_i \ot \epsilon(x_{>i}) \, y_{>i} \\ & +
		\sum (-1)^{|x|-1+|x_{<i}|} \, x_{<i}\, \epsilon(y_{<i}) \ot \bd x_i \ast y_i \ot \epsilon(x_{>i}) \, y_{>i} \\ & +
		\sum (-1)^{|x|-1+|x_{<i}|} \, x_{<i}\, \epsilon(y_{<i}) \ot x_i \ast y_i \ot \epsilon(\bd x_{>i}) \, y_{>i}.
	\end{align*}
	Since
	\[
	\epsilon(\bd x_{>i}) = 0, \quad \bd x_i \neq 0 \Leftrightarrow |x_i| = 1,
	\]
	we have
	\begin{equation} \label{e:boundary of product 2}
		\begin{split}
			\bd x \ast y & =
			\sum (-1)^{|x|-1} \, \bd x_{<i}\, \epsilon(y_{<i}) \ot x_i \ast y_i \ot \epsilon(x_{>i})\, y_{>i} \\ & +
			\sum x_{<i}\, \epsilon(y_{<i}) \ot \bd x_i \ast y_i \ot \epsilon(x_{>i})\, y_{>i}.
		\end{split}
	\end{equation}
	We also have
	\begin{align*}
		(-1)^{|x|} \, x \ast \bd y & =
		\sum x_{<i} \, \epsilon(\bd y_{<i}) \ot x_i \ast y_i \ot \epsilon(x_{>i})\, y_{>i} \\ & +
		\sum (-1)^{|y_{<i}|} \, x_{<i}\, \epsilon(y_{<i}) \ot x_i \ast \bd y_i \ot \epsilon(x_{>i}) \, y_{>i} \\ & +
		\sum (-1)^{|y_{<i}| + |y_i|} \, x_{<i}\, \epsilon(y_{<i}) \ot x_i \ast y_i \ot \epsilon(x_{>i}) \, \bd y_{>i},
	\end{align*}
	which is equivalent to
	\begin{equation} \label{e:boundary of product 3}
		\begin{split}
			(-1)^{|x|} \, x \ast \bd y & =
			\sum x_{<i} \, \epsilon(y_{<i}) \ot x_i \ast \bd y_i \ot \epsilon(x_{>i})\, y_{>i} \\ & +
			\sum x_{<i}\, \epsilon(y_{<i}) \ot x_i \ast y_i \ot \epsilon(x_{>i})\, \bd y_{>i}.
		\end{split}
	\end{equation}
	Putting identities \eqref{e:boundary of product 1}, \eqref{e:boundary of product 2} and \eqref{e:boundary of product 3} together, we get
	\begin{multline*}
		\bd (x \ot y) \ +\ \bd x \ast y\ +\, (-1)^{|x|}x \ast \bd y \\
		= \sum \epsilon(y_{<i})\, x_{<i} \ot \big(\bd(x_i \ast y_i) + \bd x_i \ast y_i + x_i \ast \bd y_i\big) \ot \epsilon(x_{>i})\, y_{>i}.
	\end{multline*}
	Since
	\begin{align*}
		\bd(x_i \ast y_i)\ +\ \bd x_i \ast y_i\ +\ x_i \ast \bd y_i =
		\epsilon(x_i)y_i\ -\ \epsilon(y_i)x_i,
	\end{align*}
	we have
	\begin{multline*}
		\bd (x \ast y) \ +\ \bd x \ast y\ +\ (-1)^{|x|}x \ast \bd y = \\
		\sum \epsilon(y_{<i}) \, x_{<i} \ot \epsilon(x_{\geq i}) y_{\geq i}\ -\
		\epsilon(y_{\leq i}) \, x_{\leq i} \ot \epsilon(x_{>i}) y_{>i} \\ =
		\epsilon(x)y - \epsilon(y)x,
	\end{multline*}
	as desired, where the last equality follows from a telescopic sum argument.
\end{proof}

\subsection{$\texorpdfstring{E_\infty}{E-infty}$-coalgebra on cubical chains}

\cref{l:cubical chain bialgebra} defines a functor from the cube category to that of $\M$-bialgebras.
This category is not cocomplete so we do not expect to have an $\M$-bialgebra structure on arbitrary cubical sets.
For example, consider the chains on the cubical set $X$ whose only non-degenerate simplices are $v, w \in X_0$.
By degree reasons $v \ast w = 0$ for any degree 1 product $\ast$ in $\chains(X)$.
The third relation in $\M$ would then imply the contradiction $0 = w-v$.
Since categories of coalgebras over operads are cocomplete we have the following.

\begin{theorem} \label{t:lift to e infinity coalgebras}
	The Yoneda extension of the composition of the functor $\cube \to \biAlg_{\M}$ defined in \cref{l:cubical chain bialgebra} with the forgetful functor $\biAlg_{\M} \to \coAlg_{\UM}$ endows the chains of a cubical set with a natural $E_\infty$-coalgebra extension of the Serre coalgebra structure.
\end{theorem}

\subsection{Cohomology operations}

In \cite{steenrod1947products}, Steenrod introduced natural operations on the mod~2 cohomology of spaces, the celebrated \textit{Steenrod squares}
\[
\begin{tikzcd} [column sep=small, row sep=0]
	\Sq^k \colon &[-20] \rH^{-n} \arrow[r] & \rH^{-n-k} \\ &
	{[\alpha]} \arrow[r, mapsto] & \big[ (\alpha \ot \alpha) \Delta_{n-k} \big],
\end{tikzcd}
\]
via an explicit construction of natural linear maps $\Delta_i \colon \chains(X) \to \chains(X) \ot \chains(X)$ for any simplicial set $X$, satisfying up to signs the following homological relations
\begin{equation} \label{e:cupi homological relations}
\bd \circ \, \Delta_i + \Delta_i \circ \bd =
(1 + T) \Delta_{i-1},
\end{equation}
with the convention $\Delta_{-1} = 0$.
These so-called \textit{cup-$i$ coproducts} appear to be fundamental.
We mention two results supporting this claim.
In higher category theory they define the nerve of $n$-categories \cite{medina2020globular} as introduced by Street \cite{street1987orientals}; and, in connection with K- and L-theory, the Ranicki--Weiss assembly \cite{ranicki1990assembly} can be used to show that chain complex valued presheaves over a simplicial complex $X$ can be fully faithfully modeled by comodules over the symmetric coalgebra structure they define on $\chains(X)$ \cite{medina2022assembly}.

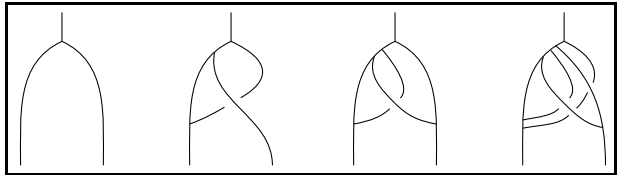
\begin{figure}[h!]
	\centering
	\boxed{
	\begin{tikzpicture}[scale=.55]
		\draw (1,3.7) to (1,3);
		\draw (1,3) to [out=205, in=90] (0,0);
		\draw (1,3) to [out=-25, in=90] (2,0);
	\end{tikzpicture}
	\hspace*{1cm}
	\begin{tikzpicture}[scale=.55]
		\draw (1,3.7) to (1,3);
		\draw (1,3) to [out=205, in=90] (0,0);
		\draw [shorten >= 0cm] (.6,2.73) to [out=-100, in=90] (2,0);
		\draw [shorten >= .15cm] (1,3) to [out=-25, in=30, distance=1.1cm] (1,1.5);
		\draw [shorten <= .1cm] (1,1.5) to [out=210, in=20] (0,1);
	\end{tikzpicture}
	\hspace*{1cm}
	\begin{tikzpicture}[scale=.55]
		\draw (1,3.7) to (1,3);
		\draw (1,3) to [out=205, in=90] (0,0);
		\draw (1,3) to [out=-25, in=90] (2,0);
		\draw [shorten >= 0cm] (.5,2.63) to [out=-110, in=135] (1,1.5);
		\draw [shorten <= 0cm] (1,1.5) to [out=-45, in=170] (2,1);
		\draw [shorten >= .1cm] (.69,2.8) to [out=-50, in=45] (1,1.5);
		\draw [shorten <= .1cm] (1,1.5) to [out=-135, in=10] (0,1);
	\end{tikzpicture}
	\hspace*{1cm}
	\begin{tikzpicture}[scale=.55]
		\draw (1,3.7) to (1,3);
		\draw (1,3) to [out=205, in=90] (0,0);
		\draw [shorten >= 0cm] (.5,2.63) to [out=-110, in=135] (1,1.5);
		\draw [shorten >= .05cm] (1,1.5) to [out=-45, in=170] (2.02,.9);
		\draw [shorten >= .1cm] (.68,2.78) to [out=-50, in=45] (1,1.5);
		\draw [shorten <= .1cm] (1,1.5) to [out=-135, in=10] (0,1.1);
		\draw (.79,2.9) to [out=-40, in=90] (2,0);
		\draw (1,3) to [out=-25, in=70] (1.7,2);
		\draw [shorten <= .15cm] (1.7,2) to [out=-120, in=45] (1.3,1.38);
		\draw [shorten <= .1cm] (1.24,1.34) to [out=-135, in=10] (0,.9);
	\end{tikzpicture}
}
	\caption{Graphs representing cup-$i$ coproducts.}
	\label{f:cup-i}
\end{figure}

In the cubical case, cup-$i$ coproducts were defined in \cite{kadeishvili1999coproducts} and \cite{pilarczyk2016cubical}.
The formulas used by these authors are similar to those introduced in \cite{medina2022fast_sq} for the simplicial case, a dual yet equivalent version of Steenrod's original.
A new description of cubical cup-$i$ coproducts can be deduced from our $E_\infty$-structure.
We first present it in a recursive form
\begin{equation} \label{e:recursive cup-i}
	\begin{split}
		& \Delta_0 = \Delta, \\
		& \Delta_i =
		(\ast \ot \id) \circ (23)(\Delta_{i-1} \ot \id) \circ \Delta.
	\end{split}
\end{equation}
A closed form formula for $\Delta_i$ uses the $\big(\ceil*{\frac{i+2}{2}}, \floor*{\frac{i+2}{2}}\big)$-shuffle permutation $\sigma_{i+2} \in \sym_{i+2}$ mapping the first and second ``decks'' to odd and even integers respectively.
Explicitly, this shuffle permutation is defined by
\begin{equation} \label{e:shuffle perm}
	\sigma_{i+2}(\ell) =
	\begin{cases}
		2\ell-1 & \ell \leq \ceil*{\frac{i+2}{2}}, \\
		2(\ell-\ceil*{\frac{i+2}{2}}) & \ell > \ceil*{\frac{i+2}{2}}.
	\end{cases}
\end{equation}
Let $\Delta^0 = \ast^0 = \id$ and define for any $k \in \N$
\begin{equation} \label{e:iterated comb}
	\begin{split}
		\ast^{k+1} &= \ast \circ (\ast^k \ot \id), \\
		\Delta^{k+1} &= (\Delta^k \ot \id) \circ \Delta.
	\end{split}
\end{equation}
With this notation it can be checked that \cref{e:recursive cup-i} is equivalent to
\begin{equation} \label{e:closed cup-i}
	\Delta_i = \left(\ast^{\ceil*{\frac{i+2}{2}}} \ot \ast^{\floor*{\frac{i+2}{2}}}\right) \circ \sigma_{i+2}^{-1} \, \Delta^{i+1}.
\end{equation}
The first four cup-$i$ coproducts are the images in the endomorphism operad of cubical (and simplicial) chains of the elements $\UM$ represented by the graphs in \cref{f:cup-i}.

It is not known if the cup-$i$ coproducts defined in \cref{e:closed cup-i} agree with those previously constructed, for which a comparison is also missing.
This highlights the value of a potential axiomatic characterization of cubical cup-$i$ coproducts as it exists in the simplicial case \cite{medina2022axiomatic}.

As already mentioned, cup-$i$ coproducts represent the Steenrod squares at the chain level,
which are primary operations in mod~2 cohomology.
To obtain secondary cohomology operations one studies the cohomological relations these operations satisfy, for example the Cartan and Adem relations \cite{steenrod1962cohomology}.
To do this at the cubical cochain level, as it was done in \cite{medina2020cartan,medina2021adem} for the simplicial case, the operadic viewpoint is important, so our $E_\infty$-structure on cubical cochains invites the construction of cochain representatives for secondary operations in the cubical case.

For $p$ an odd prime, Steenrod also introduced operations on the mod~$p$ cohomology of spaces using the homology of symmetric groups \cite{steenrod1952reduced, steenrod1953cyclic}.
Using the operadic framework of May \cite{may1970general}, we described in \cite{medina2021may_st} elements in $\UM$ representing multicooperations defining Steenrod operations at any prime.
In particular, as proven in this work, these so-called \textit{cup-$(p,i)$ coproducts} are defined on cubical chains and are expressible, similarly to \cref{e:closed cup-i}, in terms of $\Delta$, the permutations of factors, and $\ast$.
The aforementioned construction of cubical cup-$(p,i)$ coproducts has been implemented in the open source computer algebra system \href{https://comch.readthedocs.io/en/latest/}{\texttt{ComCH}} \cite{medina2021comch}.

	\section{The Cartan--Serre map} \label{s:comparison}

Let us consider, with their usual CW structures, the topological simplex $\gsimplex^n$ and the topological cube $\gcube^n$.
In \cite[p. 442]{serre1951homologie}, Serre described a quasi-isomorphism of coalgebras between the simplicial and cubical singular chains of a topological space.
It is given by precomposing with a canonical cellular map $\cs \colon \gcube^n \to \gsimplex^n$ also considered in \cite[p.199]{eilenberg1953acyclic} where it is attributed to Cartan.

The goal of this section is to deduce from a more general categorical statement that this comparison map between singular chains of a space is a quasi-isomorphism of $E_\infty$-coalgebras.

\subsection{Simplicial sets} \label{ss:simplicial sets}

We denote the \textit{simplex category} by $\simplex$, the category $\Fun(\simplex^\op, \Set)$ of \textit{simplicial sets} by $\sSet$, and the representable simplicial set $\yoneda \big( [n] \big)$ by $\simplex^n$.
As usual, we denote an element in $\simplex^n_m$ by a non-decreasing tuple $[v_0, \dots, v_m]$ with $v_i \in \{0, \dots, n\}$.
The \textit{Cartesian product} of simplicial sets is defined by the product of functors.
The \textit{simplicial $n$-cube} $\scube{n}$ is the $n^\th$-fold Cartesian product of $\simplex^1$ with itself.

We will use the following model of the topological $n$-simplex:
\[
\gsimplex^n = \big\{ (y_1, \dots, y_n) \in \gcube^n \mid i \leq j \Rightarrow y_i \geq y_j \big\},
\]
whose cell structure associates $[v_0, \dots, v_m]$ with the subset
\begin{equation} \label{e:cell structure of gsimplex}
	\Big\{ \big( \underbrace{1, \dots, 1}_{v_0}, \underbrace{y_1, \dots y_1}_{v_1-v_0}, \dots, \underbrace{y_m, \dots y_m}_{v_m-v_{m-1}}, \underbrace{0, \dots, 0}_{n-v_m} \big) \mid y_1 \geq \dots \geq y_m \Big\}.
\end{equation}
The spaces $\gsimplex^n$ define a functor $\simplex \to \CW$ with
\begin{align*}
	\sigma_i(x_1, \dots, x_n) &= (x_1, \dots, \widehat x_i, \dots, x_n) \\
	\delta_0(x_1, \dots, x_n) &= (1, x_1, \dots, x_n), \\
	\delta_i(x_1, \dots, x_n) &= (x_1, \dots, x_i, x_i, \dots, x_n), \\
	\delta_n(x_1, \dots, x_n) &= (x_1, \dots, x_n, 0).
\end{align*}
Its Yoneda extension is the \textit{geometric realization} functor.
It has a right adjoint $\sSing \colon \Top \to \sSet$ referred to as the \textit{simplicial singular complex} satisfying
\[
\sSing(\fZ)_n = \Top(\gsimplex^n, \fZ)
\]
for any topological space $\fZ$.

The functor of (\textit{normalized}) \textit{chains} $\schains \colon \sSet \to \Ch$ is the composition of the geometric realization functor and that of cellular chains.
We denote the composition $\schains \circ \sSing$ by $\sSchains$ and omit the superscript $\simplex$ if no confusion may result from doing so.
For any $n \in \N$, the \textit{Alexander--Whitney coalgebra structure} on $\chains(\simplex^n)$ is given by
\[
\Delta \big( [v_0, \dots, v_m] \big) = \sum_{i=0}^m [v_0, \dots, v_i] \ot [v_i, \dots, v_m],
\]
and
\[
\epsilon \big( [v_0, \dots, v_m] \big) =
\begin{cases}
	1 & \text{ if } m = 0, \\ 0 & \text{ if } m>0.
\end{cases}
\]
The degree 1 product $\ast \colon \chains(\simplex^n)^{\ot 2} \to \chains(\simplex^n)$ is defined by
\begin{align*}
	\begin{autobreak}
		\left[v_0, \dots, v_p \right] \ast
		\left[v_{p+1}, \dots, v_m\right] =
		{\begin{cases}
			(-1)^{p+|\sigma|} \left[v_{\sigma(0)}, \dots, v_{\sigma(m)}\right] &
			\text{ if } v_i \neq v_j \text{ for } i \neq j, \\
			0 & \text{ if not},
		\end{cases}}
	\end{autobreak}
\end{align*}
where $\sigma$ is the permutation that orders the totally ordered set of vertices and $(-1)^{|\sigma|}$ is its sign.
As shown in \cite[Theorem 4.2]{medina2020prop1} the assignment
\[
\counit \mapsto \epsilon, \quad \coproduct \mapsto \Delta, \quad \product \mapsto \ast,
\]
defines a natural $\M$-bialgebra on the chains of representable simplicial sets, and, by forgetting structure, also a natural $\UM$-coalgebra.
For any simplicial set, a natural $\UM$-coalgebra structure on its chains is defined by a Yoneda extension.

\subsection{The Eilenberg--Zilber maps} \label{ss:ez}

For any permutation $\sigma \in \sym_n$ let
\[
\gi_\sigma \colon \gsimplex^n \to \gcube^n
\]
be the inclusion defined by $(x_1, \dots, x_n) \mapsto (x_{\sigma(1)}, \dots, x_{\sigma(n)})$.
If $e$ is the identity permutation, we denote $\mathfrak{i}_{e}$ simply as $\mathfrak{i}$.
The maps $\set{\mathfrak{i}_\sigma}_{\sigma \in \sym_n}$ define a subdivision of $\gcube^n$ making it isomorphic to $\bars[\big]{\scube{n}}$ in $\CW$.
Using this identification, the identity map induces a cellular map
\[
\ez \colon \gcube^n \to \bars[\big]{\scube{n}}.
\]
We denote the induced chain map by
\[
\EZ \colon \chains(\cube^n) \to \chains\big(\scube{n}\big).
\]
For any topological space $\fZ$, the cubical map
\[
\cubify\sSing(\fZ) \to \cSing(\fZ)
\]
is defined, using the adjunction isomorphism
\[
\sSet\big(\scube{n}, \sSing(\fZ)\big) \cong
\Top\big(\bars{\scube{n}}, \fZ\big),
\]
by the assignment
\[
\big( \bars{\scube{n}} \xra{f} \fZ \big) \mapsto
\big( \gcube^n \xra{\ez} \bars{\scube{n}} \xra{f} \fZ \big).
\]
We denote the induced chain map by
\[
\EZ_{\Schains(\fZ)} \colon
\cchains\!\big(\cubify\sSing(\fZ)\big) \to \cSchains(\fZ).
\]



\subsection{The Cartan--Serre maps}

The cellular map
\[
\cs \colon \gcube^n \to \gsimplex^n
\]
is defined by
\[
\cs(x_1, \dots, x_n) =
(x_1,\ x_1 x_2, \, \dots \, , \ x_1 x_2 \dotsm x_n).
\]
We denote its induced chain map by
\[
\CS \colon \chains(\cube^n) \to \chains(\simplex^n).
\]
The chain map
\[
\CS_{\Schains(\fZ)} \colon \sSchains(\fZ) \to \cSchains(\fZ)
\]
between the singular chain complexes of a topological space $\fZ$ is defined by
\[
\CS_{\Schains(\fZ)} (\gsimplex^n \to \fZ) =
(\gcube^n \xra{\cs} \gsimplex^n \to \fZ).
\]

These maps were considered in \cite[p. 442]{serre1951homologie} where it was stated that $\CS_{\Schains(\fZ)}$ is a natural quasi-isomorphisms of coalgebras.
We will prove this in \cref{ss:e-infty preservation} showing in fact that it is a quasi-isomorphism of $E_\infty$-coalgebras.

\subsection{No-go results}

Since $\CS$ is shown to be a coalgebra map in \cref{l:cs coalgebra map} and $\EZ$ is well known to be one, one may hope for higher structures to be preserved by these maps.
We now provide some examples constraining the scope of these expectations.

\begin{example*}
	We will show that $\EZ$ does not preserve $\UM$-structures.
	More specifically, that in general
	\[
	\EZ^{\ot 2} \circ \, \Delta_1 \neq \Delta_1 \circ \EZ
	\]
	where
	\[
	\Delta_1 = (\ast \ot \id) \circ (\id \ot (12) \Delta) \circ \Delta
	\]
	is the cup-1 coproduct presented in \cref{e:closed cup-i}.
	Up to signs, on one hand we have
	\begin{align*}
		\begin{autobreak}
			\Delta_1 \big( [01] [01] \big)
			= [01][01] \ot [1][01]
			\ +\ [01][1] \ot [01][01]
			\ +\ [0][01] \ot [01][01]
			\ +\ [01][01] \ot [01][0].
		\end{autobreak}
	\end{align*}
	Therefore,
	\begin{align*}
		\begin{autobreak}
			\EZ^{\ot 2} \circ \, \Delta_1 \big( [01][01] \big)
			= \big( 011\times001 + 001\times011 \big) \ot 11 \times 01
			\ +\ 01\times11 \ot \big( 011\times001+ 001\times011 \big)
			\ +\ 00\times01 \ot \big( 011\times001 + 001\times011 \big)
			\ +\ \big( 011\times001 + 001\times011 \big) \ot 01\times00.
		\end{autobreak}
	\end{align*}
	On the other hand, we have
	\[
	\Delta_1 [0,1,2] = [0,1,2] \ot [0,1] + [0,2] \ot [0,1,2] + [0,1,2] \ot [1,2].
	\]
	Therefore,
	\begin{align*}
		\begin{autobreak}
			\Delta_1 \circ \EZ \big( [01] [01] \big)
			= \Delta_1 \big( 011\times001 + 001\times011 \big)
			= 011\times001 \ot 01\times00
			\ + \ 01\times01 \ot 011\times001
			\ + \ 011\times001 \ot 11\times01
			\ + \ 001\times011 \ot 00\times01
			\ + \ 01\times01 \ot 001\times011
			\ + \ 001\times011 \ot 01\times11.
		\end{autobreak}
	\end{align*}
	We conclude that
	\[
	\EZ^{\ot 2} \circ \, \Delta_1 \big( [01][01] \big) \neq \Delta_1 \circ \EZ \big( [01] [01] \big)
	\]
	since, for example, the basis element $01\times11 \ot 011\times001$ appears in the left sum but not in the right one.
\end{example*}

\begin{example*}
	We will show that the Cartan--Serre map does not preserve $\M$-structures.
	More specifically, that in general
	\[
	\CS(x \ast y) \neq \CS(x) \ast \CS(y).
	\]
	Consider $x = [1][1]$ and $y = [0][01]$.
	On one hand we have that
	\[
	\CS \big( [1][1] \big) \ast \CS \big( [0][01] \big) = 0
	\]
	since $\CS \big( [0][01] \big) = 0$.
	On the other hand we have, up to a signs, that
	\[
	\CS \Big( ([1][1]) \ast ([0][01]) \Big) =
	\CS \Big( [01][01] \Big) = [012],
	\]
	which establishes the claim.
\end{example*}

The reason for this incompatibility is that $\ast$ in the simplicial context is commutative, which is not the case in the cubical one.

\begin{example*}
	We will show that the Cartan--Serre map does not preserve $\UM$-structures.
	More specifically, that in general
	\[
	\CS \circ\, \widetilde\Delta_1 \neq
	\widetilde\Delta \circ \CS
	\]
	where
	\[
	\widetilde\Delta_1 =
	(\ast \ot \id) \circ (12) (\id \ot (12) \Delta) \circ \Delta.
	\]
	On one hand we have that
	\[
	\CS\Big( \widetilde\Delta_1 \big( [01][01] \big) \Big) =
	T \Delta_1 \big( [012] \big),
	\]
	and on the other that
	\[
	\widetilde\Delta_1 \circ \CS \big( [01][01] \big) =
	\Delta_1 \big( [012] \big),
	\]
	which establishes the claim.
\end{example*}

In \cref{ss:e-infty preservation} we will show that $\CS$ is a morphism of $E_\infty$-coalgebras.
To do so we now introduce an $E_\infty$-suboperad of $\UM$ where the incompatibility resulting from the lack of commutativity of $\ast$ in the cubical context is dealt with.

\subsection{Shuffle graphs} \label{ss:shuffle graphs}

Consider $k = k_1+\dots+k_r$.
A $(k_1,\dots,k_r)$-shuffle $\sigma$ is a permutation in $\sym_{k}$ satisfying
\begin{align*}
	&\sigma(1) < \dots < \sigma(k_1), \\
	&\sigma(k_1+1) < \dots < \sigma(k_1+k_2), \\
	&\qquad \vdots \\
	&\sigma(k-k_r+1) < \dots < \sigma(k).
\end{align*}
The (\textit{left comb}) \textit{shuffle graph} associated to such $\sigma$ is the $(1,k)$-graph
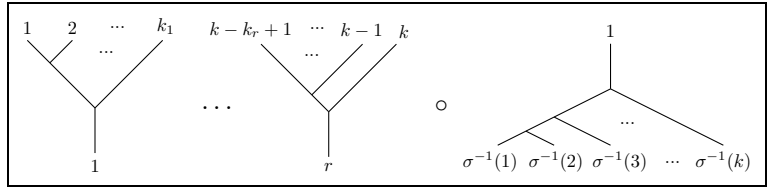
\[
\boxed{
	\begin{tikzpicture}[scale=1.5]
		\coordinate (o) at (0,0);
		\draw (o)--(0,-.4) node[below, scale=.7]{$1$};
		\draw (o)--(-.6,.6) node[above, scale=.7]{$1$\hspace{0pt}};
		\draw (-.4,.4)--(-.2,.6) node[above, scale=.7]{$2$\hspace{0pt}};
		\draw (o)--(.6,.6) node[above, scale=.7]{\hspace{0pt} $k_1$};
		\node[scale=.7] at (.1,.5){...};
		\node[scale=.7] at (.2,.72){...};
		\node at (1.1,0){$\dots$};
	\end{tikzpicture}
	\hspace*{-16pt}
	\begin{tikzpicture}[scale=1.5]
		\coordinate (o) at (0,0);
		\draw (o)--(0,-.4) node[below, scale=.7]{$r$};
		\draw (o)--(-.6,.6) node[above, scale=.7]{$k-k_r+1$\hspace*{11pt}};
		\draw (-.15,.15)--(.3,.6) node[above, scale=.7]{$k-1$\hspace*{1pt}};
		\draw (o)--(.6,.6) node[above, scale=.7]{\hspace*{5pt} $k$};
		\node[scale=.7] at (-.15,.5){...};
		\node[scale=.7] at (-.1,.74){...};
		\node at (1,.05){$\circ$};
	\end{tikzpicture}
	\begin{tikzpicture}[scale=1.5]
		\coordinate (o) at (0,0);
		\draw (o)--(0,.4) node[scale=.7, above]{$1$};
		\draw (o)--(-1,-.5) node[scale=.7, below]{$\sigma^{-1}(1)$\hspace*{8pt}};
		\draw (-.75,-.375)--(-.5,-.5) node[scale=.7, below]{\hspace*{2pt}$\sigma^{-1}(2)$};
		\draw (-.5,-.25)--(0,-.5) node[scale=.7, below]{\hspace*{10pt}$\sigma^{-1}(3)$};
		\draw (o)--(1,-.5) node[scale=.7, below]{$\sigma^{-1}(k)$};
		\node[scale=.7] at (.15,-.3){...};
		\node[scale=.7] at (.55,-.65){...};
	\end{tikzpicture}
}
\]
presented as a composition of (left comb) self-graftings of the generators $\product$ and $\coproduct$.
With the notation introduced in \cref{e:iterated comb}, the $\UM$-coalgebra sends the shuffle graph associated to $\sigma$ to
\[
(\ast^{k_1} \ot \dotsb \ot \ast^{k_r}) \circ \sigma^{-1} \Delta^{k-1}.
\]

\begin{example*}
	All the graphs in \cref{f:cup-i} are shuffle graphs.
	In fact, all the cup-$i$ coproducts presented in \cref{e:closed cup-i} are induced from shuffle graphs, whereas
	\begin{align*}
		\widetilde\Delta_1 &=
		(\ast \ot \id) \circ (12) (\id \ot (12) \Delta) \circ \Delta \\ &=
		(\ast \ot \id) \circ (123) \Delta^2,
	\end{align*}
	used in the previous section to probe the limits of the structure preserving properties of $\CS$, is not.
\end{example*}

The operad $\UMsh$ is defined as the suboperad of $\UM$ (freely) generated by shuffle graphs.
Explicitly, any element in $\UMsh(r)$ is represented by a linear combination of $(1,r)$-graphs obtained by grafting these.
The same proof used in \cite[p.5]{medina2020prop1} to show that $\UM$ is an $E_\infty$-operad can be used to prove the same for $\UMsh$.

\subsection{$E_\infty$-coalgebra preservation} \label{ss:e-infty preservation}

We devote this subsection to the proof of the following key result.

\begin{theorem} \label{t:main local}
	The chain map $\CS \colon \chains(\cube^n) \to \chains(\simplex^n)$ is a quasi-isomorphism of $\UMsh$-coalgebras.
\end{theorem}

We start by stating an alternative description of the $\CS$ map.

\begin{lemma} \label{l:cs explicit}
	Let $x = x_1 \ot \cdots \ot x_n \in \chains(\cube^n)_m$ be a basis element with $x_{q_i} = [0,1]$ for all $\{q_1 < \dots < q_m\}$.
	If there is $x_\ell = [0]$ with $\ell < q_m$ then $\CS(x) = 0$, otherwise
	\[
	\CS(x) = \big[ q_1-1, \ \dots \, , \ q_m-1, \ p(x)-1 \big]
	\]
	where $p(x) = \min \set[\big]{\ell \mid x_\ell = [0]}$ or $p(x) = n+1$ if this set is empty.
\end{lemma}

\begin{proof}
	This can be directly verified using the cell structure of $\gsimplex^n$ described in \cref{e:cell structure of gsimplex}.
\end{proof}

\begin{lemma} \label{l:cs coalgebra map}
	The chain map $\CS \colon \chains(\cube^n) \to \chains(\simplex^n)$ is a quasi-isomorphism of coalgebras.
\end{lemma}

\begin{proof}
	The chain map $\CS$ is a quasi-isomorphism compatible with the counit since it is induced from a cellular map between contractible spaces.
	We need to show it preserves coproducts.
	By naturality it suffices to verify this on $[0,1]^{\ot n}$.
	Recall from \cref{l:coproduct description} that
	\[
	\Delta \big( [0,1]^{\ot n} \big) =
	\sum_{\lambda \in \Lambda} (-1)^{\ind \lambda} \
	\Big(x_1^{(\lambda)} \ot \cdots \ot x_n^{(\lambda)}\Big) \ot
	\Big(y_1^{(\lambda)} \ot \cdots \ot y_n^{(\lambda)}\Big),
	\]
	where the sum is over all choices for each $i \in \{1,\dots,n\}$ of
	\begin{align*}
		x_i^{(\lambda)} &= [0,1],&&\text{or} & x_i^{(\lambda)} &= [0], \\
		y_i^{(\lambda)} &= [1], && & y_i^{(\lambda)} &= [0,1].
	\end{align*}
	By \cref{l:cs explicit}, the summands above not sent to $0$ by $\CS \ot \CS$ are those basis elements for which $x_i^{(\lambda)} = [0]$ implies $x_j^{(\lambda)} = [0]$ for all $i < j$.
	For any one such summand, its sign is positive and its image by $\CS \ot \CS$ is $[0, \dots, k] \ot [k, \dots, n]$ where $k+1 = \min\set[\big]{i \mid x_i^{(\lambda)} = [0]}$ or $k = n$ if this set is empty.
	The summands $[0, \dots, k] \ot [k,\dots,n]$ are precisely those appearing when applying the Alexander--Whitney coproduct to $\CS \big([0,1]^{\ot n}\big) = [0,\dots,n]$.
	This concludes the proof.
\end{proof}

We will consider the basis of $\chains(\cube^n)$ as a poset with
\[
(x_1 \ot \cdots \ot x_n) \leq (y_1 \ot \cdots \ot y_n)
\]
if and only if $x_\ell \leq y_\ell$ for each $\ell \in \{1, \dots, n\}$ with respect to
\[
[0] < [0,1] < [1].
\]
As we prove next, an example of ordered elements are the tensor factors of each summand in the iterated Serre diagonal.

\begin{lemma} \label{l:order iterated coproduct}
	Writing
	\[
	\Delta^{k-1} \big([0,1]^{\ot n}\big) =
	\sum \pm \ x^{(1)} \ot \cdots \ot x^{(k)}
	\]
	with each $x^{(\ell)}$ a basis element of $\chains(\cube^n)$, we have
	\[
	x^{(1)} \leq \cdots \leq x^{(k)}
	\]
	for every summand.
\end{lemma}

\begin{proof}
	This can be proven using a straightforward induction argument whose base case follows from inspecting \cref{l:coproduct description}.
\end{proof}

\begin{lemma}
	Let $x$, $y$ and $z$ be basis elements of $\chains(\cube^n)$.
	If both $x \leq z$ and $y \leq z$ then either $(x \ast y) = 0$ or every summand in $(x \ast y)$ is $\leq z$.
\end{lemma}

\begin{proof}
	Recall that
	\begin{multline*}
		(x_1 \ot \cdots \ot x_n) \ast (y_1 \ot \cdots \ot y_n)
		=
		(-1)^{|x|} \sum_{\ell=1}^n x_{<\ell}\, \epsilon(y_{<\ell}) \ot x_\ell \ast y_\ell \ot \epsilon(x_{>\ell}) \, y_{>\ell}.
	\end{multline*}
	By assumption $x_{<\ell} \leq z_{<\ell}$ and $y_{>\ell} \leq z_{>\ell}$ for every $\ell \in \{1, \dots, n\}$.
	If $x_\ell \ast y_\ell \neq 0$ then $x_\ell \ast y_\ell = [0,1]$ and either $x_\ell = [1]$ or $y_\ell = [1]$ which implies $z_\ell = [1]$ as well, so $x_\ell \ast y_\ell \leq z_\ell$.
\end{proof}

\begin{lemma}
	If $x$ and $y$ are basis elements of $\chains(\cube^n)$ satisfying $x \leq y$ then
	\begin{equation} \label{e:cs collapse as algebra map}
		\CS(x \ast y) = \CS(x) \ast \CS(y).
	\end{equation}
\end{lemma}

\begin{proof}
	We present this proof in the form of three claims.
	We use \cref{l:cs explicit}, the assumption $x \leq y$, and the fact that the join of basis elements in $\chains(\simplex^n)$ sharing a vertex is $0$ without explicit mention.

	\medskip\noindent \textit{Claim 1}.
	If $\CS(x) = 0$ or $\CS(y) = 0$ then for every $i \in \{1, \dots, n\}$
	\begin{equation} \label{e:zero for join}
		\CS \big( x_{<i}\, \epsilon(y_{<i}) \ot x_i \ast y_i \ot \epsilon(x_{>i}) \, y_{>i} \big) = 0.
	\end{equation}
	Assume $\CS(x) = 0$, that is, there exists a pair $p < q$ such that $x_p = [0]$ and $x_q = [0,1]$, then \eqref{e:zero for join} holds since:
	\begin{enumerate}
		\item If $i > q$, then $x_p$ and $x_q$ are part of $x_{<i}$.
		\item If $i = q$, then $x_q \ast y_q = 0$ for any $y_q$.
		\item If $i < q$, then $\epsilon(x_{>i}) = 0$.
	\end{enumerate}
	Similarly, if there is a pair $p < q$ such that $y_p = [0]$ and $y_q = [0,1]$, then \eqref{e:zero for join} holds since:
	\begin{enumerate}
		\item If $i < p$, then $y_p$ and $y_q$ are part of $y_{>i}$.
		\item If $i = p$, then $x_i = [0]$ and $x_i \ast y_i = 0$.
		\item If $i > p$, then either $x_i \ast y_i = 0$ or $x_i \ast y_i = [0,1]$ and $x_p = [0]$.
	\end{enumerate}
	This proves the first claim and identity \eqref{e:cs collapse as algebra map} under its hypothesis.

	\medskip\noindent \textit{Claim 2}.
	If $\CS(x) \neq 0$ and $\CS(y) \neq 0$ then
	\[
	\CS(x \ast y) =
	\CS \big( x_{<p_x} \epsilon(y_{<p_x}) \ot \, x_{p_x} \! \ast y_{p_x} \ot \epsilon(x_{>p_x}) \, y_{>p_x} \big)
	\]
	if $p_x = \min \big\{ i \mid x_i = [0] \big\}$ is well-defined and $x \ast y = 0$ if not.

	\medskip Assume $p_x$ is not well-defined, i.e., $x_i \neq [0]$ for all $i \in \{1, \dots, n\}$.
	Given that $x \leq y$ we have that $[0] < x_i$ implies $x_i \ast y_i = 0$, and the claim follows in this case.

	Assume $p_x$ is well-defined.
	We will show that for all $i \in \{1,\dots,n\}$ with the possible exception of $i = p_x$ we have
	\begin{equation} \label{e:case main lemma third claim}
		\CS \big( x_{<i} \, \epsilon(y_{<i}) \ot \, x_{i} \! \ast y_{i} \ot \epsilon(x_{>i}) \, y_{>i} \big) = 0
	\end{equation}
	This follows from:
	\begin{enumerate}
		\item If $i < p_x$ and $x_i = [1]$ then $y_i = [1]$ and $x_i \ast y_i = 0$.
		\item If $i < p_x$ and $x_i = [0,1]$ then $x_i \ast y_i = 0$ for any $y_i$.
		\item If $i > p_x$ then \cref{l:cs explicit} implies the claim since $x_{p_x} = [0]$ and $x_i \ast y_i \neq 0$ iff $x_i \ast y_i = [0,1]$.
	\end{enumerate}

	\noindent \textit{Claim 3}.
	If $\CS(x) \neq 0$ and $\CS(y) \neq 0$ then \eqref{e:cs collapse as algebra map} holds.

	\medskip Let us assume that $\big\{ i \mid x_i = [0] \big\}$ is empty, which implies the analogous statement for $y$ since $x \leq y$.
	Since neither of $x$ nor $y$ have a factor $[0]$ in them, \cref{l:cs explicit} implies that the vertex $[n]$ is in both $\CS(x)$ and $\CS(y)$, which implies $\CS(x) \ast \CS(y) = 0$ as claimed.

	Assume now that $p_x = \big\{ i \mid x_i = [0] \big\}$ is well defined, and let $\{q_1 < \dots < q_m\}$ with $x_{q_i} = [0,1]$ for $i \in \{1,\dots,m\}$.
	Since $\CS(x) \neq 0$ \cref{l:cs explicit} implies that $p_x > q_m$, so $\epsilon(x_{>p_x}) = 1$ and Claim 2 implies
	\[
	\CS(x \ast y) =
	\CS \big( x_{<p_x} \epsilon(y_{<p_x}) \ot \, x_{p_x} \! \ast y_{p_x} \ot y_{>p_x} \big).
	\]
	We have the following cases:
	\begin{enumerate}
		\item If $\epsilon(y_{<p_x}) = 0$ then there is $q_i$ such that $y_{q_i} = [0,1]$ so $[q_i-1]$ is in both $\CS(x)$ and $\CS(y)$.
		\item If $\epsilon(y_{p_x}) \neq 0$ and $y_{p_x} \in \{[0], [0,1]\}$ then $x_{p_x} \ast y_{p_x} = 0$ and $[p_x-1]$ is in both $\CS(x)$ and $\CS(y)$.
		\item If $\epsilon(y_{p_x}) \neq 0$ and $y_{p_x} = [1]$ let $\{\ell_1 < \dots < \ell_k\}$ be such that $y_{\ell_j} = [0,1]$ and let $p_y > \ell_k$ be either $n+1$ or $\min\{j \mid y_j = \{0\}\}$ then
		\begin{align*}
			\CS(x \ast y) & =
			\CS \big( x_{< p_x} \ot x_{p_x} \ast y_{p_x} \ot y_{> p_y} \big) \\ & =
			[q_1-1, \dots, q_m-1, p_x-1, \ell_1-1, \dots, \ell_k-1, p_y-1] \\ & =
			\CS(x) \ast \CS(y).
		\end{align*}
	\end{enumerate}
	This concludes the proof.
\end{proof}

Combining the previous two lemmas we obtain the following.

\begin{lemma} \label{l:cs product order}
	Let $x^{(1)} \leq \dots \leq x^{(k)}$ be basis elements of $\chains(\cube^n)$.
	Then,
	\[
	\CS \, \circ \ast^{k-1} \big( x^{(1)} \ot \dotsb \ot x^{(k)} \big)
	=
	\ast^{k-1} \circ \CS^{\ot k} \big( x^{(1)} \ot \dotsb \ot x^{(k)} \big).
	\]
\end{lemma}

We are now ready to present the argument establishing that $\CS$ is an $E_\infty$-coalgebra map.

\begin{proof}[Proof of \cref{t:main local}]
	Since $\UMsh$ is generated by elements represented by shuffle graphs, we only need to show that for any $(k_1,\dots,k_r)$-shuffle $\sigma$ with $k = k_1+\dots+k_r$ the following holds
	\[
	\CS^{\ot r}(\ast^{k_1} \ot \dotsb \ot \ast^{k_r}) \circ \sigma^{-1} \Delta^{k-1} =
	(\ast^{k_1} \ot \dotsb \ot \ast^{k_r}) \circ \sigma^{-1} \Delta^{k-1} \circ \CS.
	\]
	By naturality, it suffices to prove this identity for $[0,1]^{\ot n}$.
	According to \cref{l:order iterated coproduct}
	\[
	x^{(1)} \leq \cdots \leq x^{(k)}
	\]
	for every summand in
	\[
	\Delta^{k-1} \big([0,1]^{\ot n}\big) =
	\sum \pm \ x^{(1)} \ot \cdots \ot x^{(k)}.
	\]
	Since $\sigma$ is a shuffle permutation, \cref{l:cs product order} implies that
	\begin{multline*}
		\CS^{\ot r}(\ast^{k_1} \ot \dotsb \ot \ast^{k_r}) \circ \sigma^{-1} \Delta^{k-1} \big([0,1]^{\ot n}\big) \\ =
		(\ast^{k_1} \ot \dotsb \ot \ast^{k_r}) \circ \sigma^{-1} \CS^{\ot k} \circ \, \Delta^{k-1} \big([0,1]^{\ot n}\big).
	\end{multline*}
	As proven in \cref{l:cs coalgebra map}, $\CS$ is a coalgebra map, which concludes the proof.
\end{proof}

\subsection{Categorical reformulation} \label{ss:e infty preservation}

The assignment $2^n \mapsto \scube{n}$ defines a functor $\cube \to \sSet$ with
\begin{align*}
	&\delta_i^\varepsilon \colon \scube{n} \to \scube{(n+1)} \\
	&\sigma_i \colon \scube{(n+1)} \to \scube{n}
\end{align*}
given by inserting $[\varepsilon, \dots, \varepsilon]$ as the $i^\th$ factor and removing the $i^\th$ factor respectively.
Its Yoneda extension, referred to as triangulation functor, is denoted by
\[
\triangulate \colon \cSet \to \sSet.
\]
This functor admits a right adjoint
\[
\cubify \colon \sSet \to \cSet
\]
defined, as usual, by the expression
\[
\cubify(Y)(2^n) = \sSet \big( \scube{n}, \, Y \big).
\]
We mention that, as proven in \cite[\S~8.4.30]{cisinski2006presheaves}, the pair $(\triangulate,\, \cubify)$ defines a Quillen equivalence when $\sSet$ and $\cSet$ are considered as model categories.

%
%
%
%

%

\begin{definition}
	The simplicial map $\ccs \colon \scube{n} \to \simplex^n$ is defined by
	\[
	[\varepsilon_0^1, \dots, \varepsilon_m^1]
	\times \dots \times
	[\varepsilon_0^n, \dots, \varepsilon_m^n]
	\mapsto
	[v_0, \dots, v_m]
	\]
	where $v_i = \varepsilon_i^1 + \varepsilon_i^1 \varepsilon_i^2 + \dots + \varepsilon_i^1 \dotsm \varepsilon_i^n$.
\end{definition}

Please observe that the maps $\cs$ and $\bars{\ccs} \circ \ez$ agree.

\begin{definition}
	Let $Y$ be a simplicial set.
	The map
	\[
	\CS_Y \colon \schains(Y) \to \cchains(\cubify Y)
	\]
	is the linear map induced by sending a simplex $y \in Y_n$ to the composition
	\[
	\scube{n} \xra{\ccs} \simplex^n \xra{\xi_y} Y
	\]
	where $\xi_y \colon \simplex^n \to Y$ is the simplicial map determined by $\xi_y \big( [n] \big) = y$.
\end{definition}

\begin{theorem} \label{t:main comparison}
	For any simplicial set $Y$ the map $\CS_Y \colon \schains(Y) \to \cchains(\cubify Y)$ is a quasi-isomorphism of $\UMsh$-coalgebras which extend respectively the Alexander--Whitney and Serre coalgebra structures.
\end{theorem}

\begin{proof}
	This is a direct consequence of \cref{t:main local} following from a standard category theory argument, which we now present.
	Consider the isomorphism
	\[
	\chains(\cubify Y) \cong \displaystyle\bigoplus_{n\in\N} \chains(\cube^n) \ot \k\set[\Big]{\sSet \big( \scube{n}, \simplex^n \big)} \Big/ \sim \,
	\]
	and the canonical linear inclusions:
	\[
	\begin{tikzcd} [row sep=-7pt, column sep=small]
		\chains(\cube^n) \arrow[r] &
		\displaystyle\bigoplus_{m \in \N} \Hom \big( \chains(\cube^m),\, \chains(\cube^n) \big) \\
		(2^m \xra{\delta} 2^n) \arrow[r, mapsto] &
		\big(\chains(\cube^m) \xra{\chains(\delta)} \chains(\cube^n)\big)
	\end{tikzcd}
	\vspace*{-7pt}
	\]
	and
	\[
	\begin{tikzcd} [row sep=-5pt, column sep=small]
		\displaystyle \bigoplus_{n \in \N}
		\k\set[\Big]{\sSet\big(\scube{n}, \simplex^n\big)} \arrow[r] &
		\displaystyle \bigoplus_{n \in \N}
		\Hom \Big( \chains \big( \scube{n} \big),\, \chains(\simplex^n) \Big) \\
		\big( \scube{n} \xra{f} \simplex^n \big) \arrow[r, mapsto] &
		\Big( \chains \big( \scube{n} \big) \xra{\chains(f)} \chains(\simplex^n) \Big).
	\end{tikzcd}
	\]
	We can use these and the naturality of $\EZ$ to construct the following chain map which is an isomorphism onto its image.
	\[
	\begin{tikzcd}[column sep=small, row sep=0]
		\chains(\cubify Y) \arrow[r] &
		\displaystyle\bigoplus_{n\in\N} \Hom \big( \chains(\cube^n),\, \chains(Y)\, \big) \\
		(\delta \ot f) \arrow[r, mapsto]& \big( \chains(f) \circ \EZ \circ \chains(\delta) \, \big).
	\end{tikzcd}
	\]
	Let $\Gamma$ be an element in $\UMsh(r)$ and denote by $\Gamma^\cube \colon \chains(\cubify Y) \to \chains(\cubify Y)^{\ot r}$ and $\Gamma^\simplex \colon \chains(Y) \to \chains(Y)^{\ot r}$ its image in the respective endomorphism operads.
	Using the naturality of $\Gamma^\square$, we have that
	$\Gamma^\square(\delta \ot \!f)$ corresponds to $(\chains(f) \circ \EZ)^{\ot r} \circ \Gamma^\square \circ \chains(\delta)$.
	On the other hand, the map $\CS_Y$ corresponds to
	\[
	\begin{tikzcd}[column sep=small, row sep=0]
		\chains(Y)_n \arrow[r] & \chains(\cubify Y)_n \\
		y \arrow[r,mapsto]& \big( \chains(\xi_y) \circ \CS \big)
	\end{tikzcd}
	\]
	where $\xi_y \colon \simplex^n \to Y$ is determined by $\xi_y \big( [n] \big) = y$, and we used that $\CS = \chains(\ccs) \circ \EZ$ to ensure the above assignment is well defined.
	The image of $\Gamma^\triangle(y)$ corresponds to $\chains(\xi_y)^{\ot r} \circ \Gamma^\triangle \circ \CS$.
	So the claim follows from the identity
	\begin{align*}
		\Gamma^\cube \big( 2^n \ot (\xi_y \circ \ccs) \big) &=
		(\chains(\xi_y \circ \ccs) \circ \EZ)^{\ot r} \circ \Gamma^\square \\ &=
		\chains(\xi_y)^{\ot r} \circ \CS^{\ot r} \circ \, \Gamma^\square \\ &=
		\chains(\xi_y)^{\ot r} \circ \, \Gamma^\triangle \circ \CS
	\end{align*}
	where we used that $\CS^{\ot r} \circ \, \Gamma^\square = \Gamma^\triangle \circ \CS$ as proven in \cref{t:main local}.
\end{proof}

\begin{corollary} \label{c:zig-zag}
	For any cubical set $X$
	\[
	\cchains(X) \xra{\cchains(\xi_X)}
	\cchains(\cubify{\triangulate X})
	\xla{\CS_{\cT X}} \schains(\triangulate X),
	\]
	where $\xi$ is the unit of adjunction, is a natural zig-zag of quasi-isomorphisms of $\UMsh$-coalgebras which extend respectively the Serre and Alexander--Whitney coalgebra structures.
\end{corollary}

\begin{proof}
	The map $\CS_{\cT X}$ is a quasi-isomorphism of $\UMsh$-coalgebras by \cref{t:main comparison}, whereas $\cchains(\xi_X)$ is also one since it is induced from a cubical map that is a weak-equivalence.
\end{proof}

\begin{corollary} \label{c:cs e infty}
	The singular simplicial and cubical chains of a topological space $\fZ$ are quasi-isomorphic as $\UMsh$-coalgebras which extend respectively the Alexander--Whitney and Serre coalgebra structures.
	More specifically, the map
	\[
	\CS_{\Schains(\fZ)} \colon \sSchains(\fZ) \to \cSchains(\fZ)
	\]
	is a quasi-isomorphism of $\UMsh$-coalgebras.
\end{corollary}

\begin{proof}
	It can be verified using that $\cs = \bars{\ccs} \circ \ez$ that this map factors as
	\[
	\CS_{\Schains(\fZ)} \colon \sSchains(\fZ) \xra{\CS_{\sSing(\fZ)}} \cchains\!\big(\cubify\sSing(\fZ)\big) \xra{\EZ_{\Schains(\fZ)}} \cSchains(\fZ)
	\]
	where the first map was proven in \cref{t:main comparison} to be a quasi-isomorphism of $\UMsh$-coalgebras, and the second, introduced in \cref{ss:ez}, is also one since it is induced from a cubical map whose geometric realization is a homeomorphism.
\end{proof}

%
%
%

	\section{Future work} \label{s:future}

In the fifties, Adams introduced in \cite{adams1956cobar} a comparison map
\[
\cobar \sSchains(\fZ,z) \to \cSchains(\loops_z \fZ)
\]
from his cobar construction on the simplicial singular chains of a pointed space $(\fZ,z)$ to the cubical singular chains on its based loop space $\loops_z \fZ$.
This comparison map is a quasi-isomorphism of algebras, which was shown by Baues \cite{baues1998hopf} to be one of bialgebras by considering Serre's cubical coproduct.
In \cite{medina2021cobar} the $E_\infty$-coalgebra structure defined here is used to generalize Baues' result, by showing that Adams' comparison map is a quasi-isomorphism of $E_\infty$-bialgebras or, more precisely, of monoids in the category of $\UM$-coalgebras.

For a closed smooth manifold $M$, in \cite{medina2021flowing} a canonical vector field was used to compare multiplicatively two models of ordinary cohomology.
On one hand, a cochain complex generated by manifolds with corners over $M$, with partially defined intersection; on the other, the cubical cochains of a cubulation of $M$ with the Serre product.
With the explicit description introduced here of an $E_\infty$-structure on cubical cochains, we expect to build on this multiplicative comparison and enhance geometric cochains \cite{medina2022foundations} with compatible representations of further derived structure.
	\sloppy
	\printbibliography
\end{document}